\documentclass[final]{elsarticle}
\usepackage[a4paper, top=2cm, bottom=2cm, left=2.5cm, right=2.5cm]{geometry}
\usepackage{bm}
\usepackage{amssymb}
\usepackage{amsmath}
\usepackage{float}
\usepackage{graphicx}    
\usepackage{subcaption}  
\usepackage{booktabs}    
\usepackage{caption}     
\usepackage{multirow}    
\usepackage{array}       
\usepackage{makecell}    
\usepackage{nomencl} 
\usepackage[numbers]{natbib}
\usepackage{algorithm}
\usepackage{algorithmicx}
\usepackage{algpseudocode}
\usepackage{longtable}
\usepackage{xcolor}
\usepackage{mathrsfs}   

\makenomenclature
\journal{Elsevier}

\begin{document}

\begin{frontmatter}



\title{{Data-driven multifidelity and multiscale topology optimization based on phasor-based evolutionary de-homogenization}}

\author[label1]{Shuzhi Xu}

\author[label2]{Yifan Guo}

\author[label1]{Hiroki Kawabe}

\author[label1]{Kentaro Yaji\corref{cor1}}
\ead{yaji@mech.eng.osaka-u.ac.jp}
\cortext[cor1]{Corresponding author.}

\affiliation[label1]{%
  organization={Department of Mechanical Engineering, The University of Osaka},
  addressline={2-1 Yamadaoka},
  city={Suita},
  state={Osaka},
  postcode={565-0871},
  country={Japan}
}

\affiliation[label2]{%
  organization={Department of Mechanical Engineering, University of Alberta},
  city={Edmonton},
  state={Alberta},
  postcode={T6G 1H9},
  country={Canada}
}

\begin{abstract}
Multiscale topology optimization is crucial for designing porous infill structures with high stiffness-to-weight ratios and excellent energy absorption. Although gradient-based methods provide a rigorous framework, they are computationally expensive and struggle to capture cross-scale sensitivities in nonlinear settings. Moreover, the resulting hierarchical geometries are often overly complex and lack macroscopically meaningful features. To overcome these issues, we propose an evolutionary de-homogenization framework that couples MultiFidelity Topology Design (MFTD) with a phasor-based de-homogenization technique. The framework translates low-dimensional geometric descriptors into manufacturable high-resolution structures through a hybrid evolutionary algorithm integrating NSGA-II selection, VAE-enabled latent space crossover, and a novel image deformation-based mutation operator. This gradient-free approach achieves efficient optimization while ensuring geometric continuity. Numerical results confirm that the method effectively balances efficiency and design flexibility, offering a scalable pathway for fabrication-aware multiscale structural optimization.

\end{abstract}




\end{frontmatter}


\section*{Nomenclature}

\begin{longtable}{p{3cm} p{12cm}}
\toprule
\textbf{Symbol} & \textbf{Description} \\
\midrule
\( \boldsymbol{\rho}_f \) & Fine-scale material density field \\
\( \boldsymbol{\rho}_m \) & Macro-scale material density field \\
\( \mathbf{X} \) & Geometric parameter field \((\mu_1, \mu_2, \theta)\) for macro elements \\
\( \mathbf{X}^{(i)} \) & Geometric parameters for the \( i \)-th macro element \\
\( \mu_1, \mu_2 \) & Widths of horizontal and vertical bars in the Rank-2 microstructure \\
\( \theta \) & In-plane rotation angle of the microstructure \\
\( \rho \) & Macroscopic material density field for distinguishing void and solid regions \\
\( \mathbf{C}^{\text{H}} \) & Homogenized elasticity tensor in local coordinates \\
\( \tilde{\mathbf{C}}^{\text{H}} \) & Rotated homogenized elasticity tensor in global coordinates \\
\( \mathbf{T}(\theta) \) & Transformation matrix for elasticity tensor rotation \\
\( C_{11}^{\text{H}}, C_{22}^{\text{H}}, C_{12}^{\text{H}}, C_{33}^{\text{H}} \) & Components of the homogenized elasticity tensor \\
\( \alpha_{ij}^{(q,r)} \) & Polynomial fitting coefficients for elasticity interpolation \\
\( \mathbf{z} \) & Reduced design variable field after PCA \\
\( \mathbf{\Phi} \) & Principal component matrix obtained from PCA \\
\( \bar{\mathbf{X}} \) & Mean geometric field for PCA processing \\
\( \texttt{geom} \) & CAD-compatible geometry output (e.g., .dxf file) \\
\( G_{\text{map}}(\cdot) \) & Geometry mapping operator from \(\mathbf{X}\) to \texttt{geom} \\
\( \mathbf{x} \) & Spatial coordinate in computational domain \\
\( \mathcal{F} = \{ \mathbf{n}_1, \mathbf{n}_2 \} \) & Orthogonal lamination directions \\
\( \mathbf{n}_j(\theta) \) & Unit vector along the \( j \)-th lamination direction \\
\( \mathcal{W} = \{ \overline{t}_{w1}, \overline{t}_{w2} \} \) & Directional thickness fields \\
\( G_e(\mathbf{x}) \) & Phasor field emitted from coarse element \( e \) \\
\( G(\mathbf{x}, \mathbf{n}_j) \) & Synthesized global phasor field along direction \( \mathbf{n}_j \) \\
\( \tilde{\phi}(\mathbf{x}, \mathbf{n}_j) \) & Demodulated phase field along lamination direction \\
\( A_e(\mathbf{x}) \) & Sampling filter for phasor emissions \\
\( \tilde{\beta} \) & Gaussian decay parameter for phasor field kernels \\
\( \omega \) & Oscillation frequency of phasor signals \\
\( \varphi_e \) & Phase offset of phasor signal for element \( e \) \\
\( \Delta_e^{\mathbf{x}}(\mathbf{x}) \) & Anisotropic distance function used in phasor attenuation \\
\( \boldsymbol{\phi}(\mathbf{x}) \) & Binary structural field after de-homogenization \\
\( c \) & Threshold value for isocontour extraction (typically 0.5) \\
\( \phi_{00}, \phi_{10}, \phi_{01}, \phi_{11} \) & Scalar field values at the four corners of a grid cell \\
\( b_k \) & Binary value (0 or 1) assigned to each corner node based on thresholding \\
\( \mathbf{p}_{\text{end},1}, \mathbf{p}_{\text{end},2} \) & Endpoints of a cell edge for interpolation \\
\( \phi_1, \phi_2 \) & Scalar values at edge endpoints \\
\( \mathbf{p}_{\text{int}} \) & Interpolated point on the edge for isocontour \\
\( \mathcal{P} \) & Polyline representing the initial extracted boundary \\
\( s_i \) & Cumulative arc length up to point \( i \) \\
\( L \) & Total arc length of the boundary \\
\( \Delta s \) & Arc-length interval for uniform resampling \\
\( \tilde{\mathbf{p}}_k \) & Resampled point along the boundary \\
\( \tilde{\mathcal{P}} \) & Uniformly resampled boundary polyline \\
\( M \) & Simplification factor for resampling \\
\( \varepsilon_{\text{shape}} \) & Maximum geometric deviation between original and resampled polylines \\
\( \varepsilon_{\text{acc}} \) & Relative geometric accuracy (\%) after resampling \\
\( \mathbf{U} \) & Nodal displacement field in FEM \\
\( \mathbf{K} \) & Global stiffness matrix \\
\( \mathbf{F} \) & Global force vector \\
\( V(\mathbf{X}) \) & Volume function dependent on design variables \\
\( V_0 \) & Prescribed maximum volume constraint \\
\( \mathbf{O}(\mathbf{z}_i) \) & Multi-objective evaluation vector (volume fraction and structural performance) \\
\( G_{\mathrm{vf}} \) & Volume fraction objective function \\
\( G_{\mathrm{opt}} \) & Structural performance objective function (e.g., compliance, buckling load) \\
\( \mathcal{J} \) & Population set of design candidates \\
\( N_c \) & Number of principal components used in PCA \\
\( N_m \) & Number of macro elements \\
\( n \) & Number of design samples or individuals \\
\( \mathbf{z}_{\text{lat}} \) & Latent vector sampled from VAE distribution \\
\( N_{\text{lat}} \) & Dimension of the VAE latent space \\
\( \mathcal{L} \) & Loss function for training VAE \\
\( \mathcal{L}_{\text{recon}} \) & Reconstruction loss in VAE training \\
\( \mathcal{L}_{\text{KL}} \) & Kullback-Leibler divergence in VAE training \\
\( \mathcal{N}_{c} \) & Coarse mesh \\
\( \mathcal{N}_{i} \) & Intermediate mesh \\
\( \mathcal{N}_{f} \) & Fine mesh \\
\( \omega_c \) & Cubic interpolation kernel \\
\( r_m \) & The radius of perturbation \\
\( \alpha_{\mathrm{m}} \) & Attenuation factor \\
\( p_0 \) & Decay smoothness factor \\
\bottomrule
\end{longtable}

\section{Introduction} \label{sec:1}

Multiscale structures such as porous infill architectures—characterized by their ultra-lightweight properties, high stiffness, and superior energy absorption capabilities—have garnered significant attention in the fields of structural optimization and advanced manufacturing~\cite{sun2017topological, zhang2015bioinspired, guo2021recent,liu2018current,liu2025bioinspired}. Traditionally, multiscale topology optimization based on gradient-based methods has been regarded as an effective strategy for designing such structures, as it allows for the simultaneous optimization of material distribution at both the macro and micro scales~\cite{wu2021topology, xu2021multi, gao2019concurrent,yan2014concurrent,zhang2021novel,zhang2022comprehensive,liu2025clustering,liu2025novel}. Within this framework, macro-scale elements are interpreted as homogenized representations of microstructures, which are typically modeled using asymptotic homogenization theory~\cite{xia2015design, andreassen2014determine, dong2019149}, and analyzed numerically on structured grids. Despite its theoretical rigor, traditional multiscale topology optimization encounters substantial challenges. Chief among these are the heavy computational costs associated with resolving multiple scales concurrently~\cite{wu2021topology} and the difficulty in accurately capturing multiscale sensitivities, particularly in the presence of nonlinear physical phenomena~\cite{xia2017recent}. Furthermore, the optimized geometries resulting from such approaches often exhibit highly intricate forms, which impose significant obstacles for practical manufacturing.

Recent advancements in MultiFidelity Topology Design (MFTD)~\cite{yaji2020multifidelity} have opened new avenues for efficiently optimizing complex structures across diverse engineering domains. Then this method is improved by using data driven approaches~\cite{yamasaki2021data, yaji2022data}. By leveraging low-fidelity models to explore the design space and high-fidelity models to ensure accurate performance evaluations, MFTD significantly reduces computational costs while enabling the use of non-gradient-based optimization algorithms. This approach allows MFTD to address design problems that are often intractable for conventional gradient-based methods. MFTD has demonstrated strong versatility and effectiveness across a wide range of applications, including stress minimization~\cite{kato2025maximum}, heat exchanger design~\cite{kobayashi2019freeform}, thermal energy storage systems~\cite{luo2025data1}, turbulent natural convection cooling system~\cite{luo2025data}. Furthermore, recent developments—such as enhanced crossover strategies~\cite{yaji2024latent}, extensions to three-dimensional problems~\cite{yang2025data}, and the introduction of multi-field design variables~\cite{kawabe2025data}—have further expanded the scope and applicability of MFTD~\cite{kii2025data,yang2025enhanced}.

While MFTD offers a promising framework for efficient multiscale design exploration, its direct application to multiscale optimization still inevitably leads to the generation of complex hierarchical geometries, which presents significant manufacturing challenges. Besides, applying MFTD directly to single-scale porous infill structure design introduces additional difficulties: the massive number of design variables at the fine scale poses two critical issues, namely: (i) the limited capacity of generative models to handle high-dimensional data~\cite{yang2025data}, and (ii) the difficulty of achieving convergence with non-gradient-based optimization algorithms in large-scale design spaces~\cite{verma2021comprehensive}. Although the use of offline-trained deep learning models could theoretically mitigate these issues, this approach faces inherent obstacles similar to those currently encountered in the field of AI for Science, such as: (i) the requirement for a large number of high-quality training samples (which is often difficult to obtain in scientific applications)~\cite{montans2019data,wang2025exploration}, and (ii) the challenge of ensuring reliable generalization across diverse design tasks.

Meanwhile, an alternative strategy for addressing multiscale design challenges has been developed, known as de-homogenization~\cite{pantz2008post, groen2018homogenization, woldseth2024phasor}. These methods offer a computationally inexpensive yet practically feasible means of translating homogenized designs into coherent, single-scale, and manufacturable structures, thus providing a viable pathway to bridge the gap between multiscale optimization and actual fabrication. In addition, multiscale algorithms built upon the concept of de-homogenization have been successfully extended to problems involving fluid flow and convective heat transfer, further demonstrating the versatility and broad applicability of this strategy~\cite{dede2020inverse, zhou2022inverse, feppon2024multiscale, li2025multi}. Nevertheless, it should be noted that since de-homogenization is inherently a post-processing approach, it inevitably introduces geometric and performance deviations from the original homogenized design, which are often difficult to eliminate completely.

To address the aforementioned challenges, a new structural design framework has been proposed recently, termed \textit{evolutionary de-homogenization}, which seamlessly integrates MFTD with de-homogenization based geometric refinement~\cite{xu2025evolutionary}. The motivation behind this approach is straightforward: from a geometric modeling perspective, de-homogenization can be interpreted as a form of multifidelity modeling. It establishes a mapping from low-fidelity geometric descriptors—such as feature dimensions and orientations—derived from homogenized designs, to high-fidelity, detailed solid structures, in which the geometric information is significantly enriched. While directly applying gradient-based methods to multiscale optimization is often hindered by the difficulty of computing sensitivities across scales, the proposed MFTD-based approach effectively circumvents these challenges by adopting non-gradient optimization strategies.

While the previously proposed evolutionary de-homogenization framework demonstrated the feasibility of integrating multiscale design and geometric refinement, it still suffers from several limitations. First, it does not fully exploit the advantages of de-homogenization, as the lattice infill is restricted to a single microstructure parameter. Second, the approach involves a large number of design variables, leading to significant training costs for the generative model during each optimization iteration. Finally, the absence of a mutation mechanism makes the optimization highly sensitive to the quality of the initial solution.

Therefore, this study presents a substantially enhanced evolutionary de-homogenization framework, which addresses the aforementioned limitations of our previous approach. The framework introduces a phasor-based de-homogenization technique that enables more efficient and flexible geometric refinement compared to the original method. To further improve scalability, Principal Component Analysis (PCA)~\cite{abdi2010principal} is integrated for dimensionality reduction, thereby establishing an efficient mapping from a low-dimensional set of low-fidelity design variables to detailed high-fidelity manufacturable geometries. This mapping serves as the foundation of the multifidelity modeling process. Building on this foundation, an Evolutionary Algorithm (EA) is employed to optimize the reduced set of design variables~\cite{vikhar2016evolutionary}. Each candidate solution is adaptively decoded into a high-fidelity model and evaluated for structural performance. Optimized candidates are selected using the Non-Dominated Sorting Genetic Algorithm II (NSGA-II)~\cite{srinivas1994muiltiobjective}, while a Variational AutoEncoder (VAE)~\cite{kingma2013auto} is utilized to facilitate efficient crossover operations in the latent space. Finally, a image based processing method is developed to achieve the mutation operation, which could provide wider design space during optimization. Through successive iterations of selection, crossover, and mutation, the final optimized CAD model is obtained, accurately capturing both structural performance and manufacturable geometric features. The resulting designs can be directly converted into \textit{STereoLithography} (STL) format, facilitating their use in {Additive Manufacturing} (AM) processes.

The remainder of this paper is organized as follows: Section~\ref{sec:Multifidelity Method} details the proposed multifidelity modeling framework, including the abstraction of low-fidelity variables and their mapping to detailed full-scale structures. Section~\ref{sec:Optimization model} describes the overall optimization process, and each component will be explained in detail. Section~\ref{sec:Numerical Examples} presents several numerical case studies to validate the effectiveness of the proposed framework. Finally, discussions and conclusions are provided in Section~\ref{sec:Concluison}.

\section{Multifidelity model formulation} \label{sec:Multifidelity Method}

The core concept of Evolutionary de-homogenization lies in adopting a multifidelity modeling framework to establish an explicit, knowledge-driven mapping between a set of simple control variables and highly complex geometric structures. By adjusting these simple variables, one can induce systematic and meaningful variations in the resulting geometries without directly manipulating the fine-scale structural details.

Within this framework, the low-fidelity model refers to computationally efficient simulations performed on homogenized representations, where the design is described using a small set of control variables. In contrast, the high-fidelity model involves detailed simulations on fully resolved geometries, capturing intricate features through adaptive meshing and fine-scale resolution.

\begin{figure}[!tbp]
  \centering
  \includegraphics[width=\textwidth]{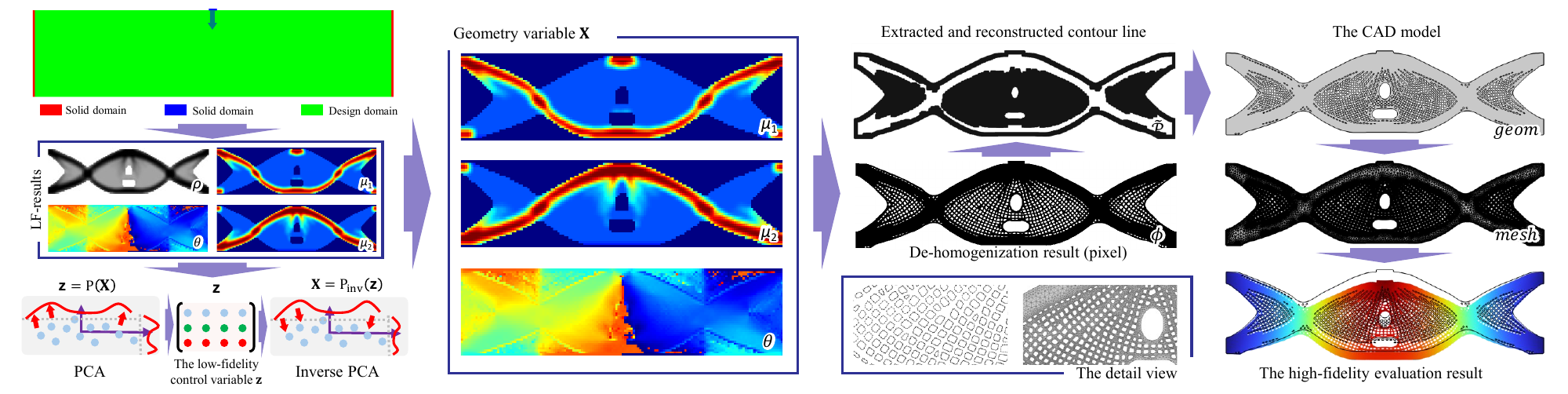}
  \caption{The whole procesure for the multifidelity modelling, in which a double clamped beam is taken as an example.}
  \label{fig:MFmodel_pipeline}
\end{figure}

The subsequent sections present a comprehensive introduction of this multi-fidelity modeling strategy and elaborate on the methodology used to construct the mapping from simple design parameters to manufacturable complex structures. As illustrated in Figure~\ref{fig:MFmodel_pipeline}, a double-clamped beam structure is used as a representative example to conceptually demonstrate this process.

\subsection{Overview of multifidelity model} \label{subsec:2.1}

The construction of the multifidelity model in this work  is based on the homogenization theory, which provides a rigorous mathematical framework for deriving effective macroscopic properties from periodic microstructures. In the context of density-based topology optimization, the material distribution is represented by a fine-scale density vector \( \boldsymbol{\rho}_\text{f} \in [0,1]^{N_\text{f}} \), where each entry \( \rho_{\text{f},e} \)  denotes the material density in the \( e \)-th finite element. While binary values (\( \rho_e = 0 \) or \( 1 \)) represent void and solid states, continuous relaxation is commonly adopted to ensure numerical tractability. According to homogenization theory, the behavior of a periodic microstructure can be equivalently described by a homogeneous material with effective properties. Thus, the fine-scale field \( \boldsymbol{\rho}_\text{f} \) can be coarse-grained into a macro-scale density vector \( \boldsymbol{\rho}_\text{m} \in [0,1]^{N_\text{m}} \), with each macro element encoding the response of a corresponding representative microstructure. This multiscale modeling strategy significantly reduces computational effort in macroscopic simulations while preserving the influence of microscale design features.

\begin{figure}[!tbp]
  \centering
  \includegraphics[width=\textwidth]{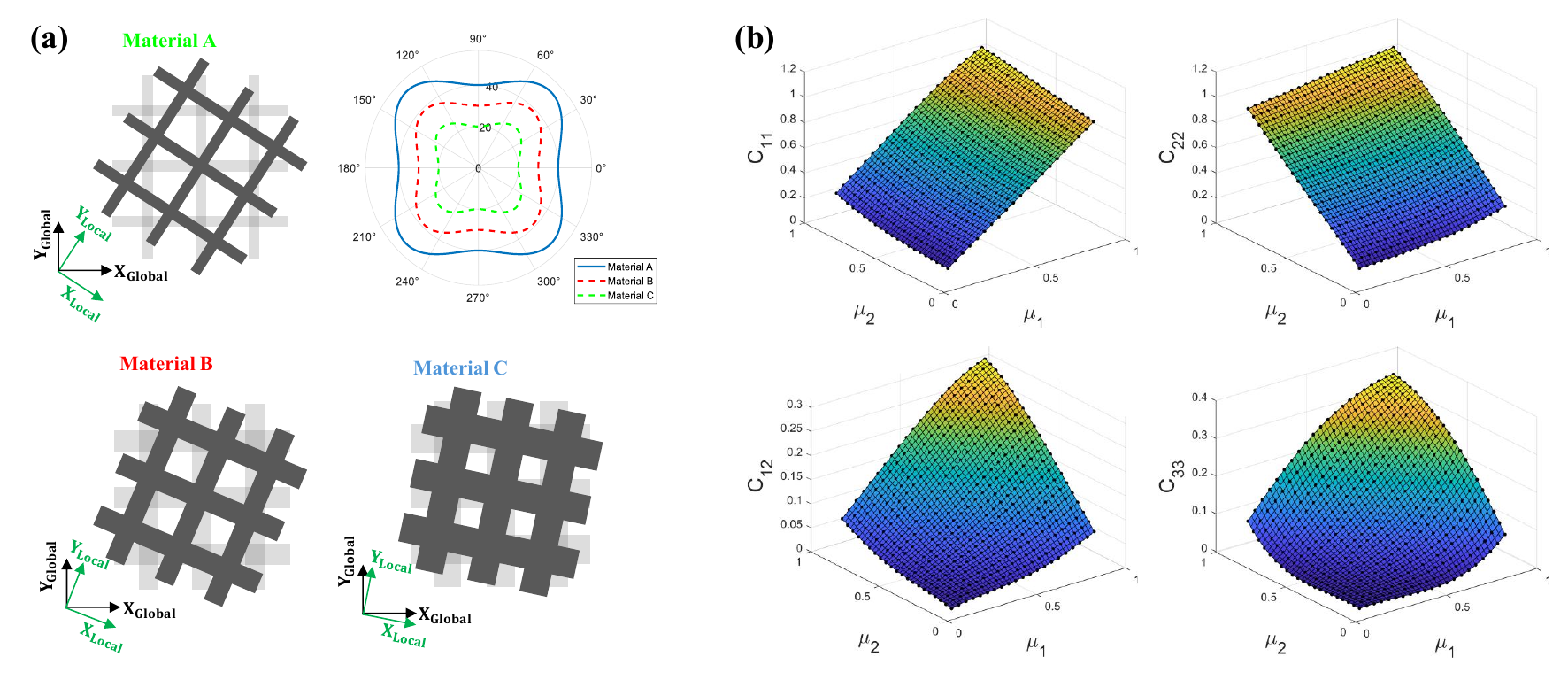}
  \caption{The illustration of the considered microstructures: (a) the orthotropic properties of composite material and its corresponding geometry configuration; (b) the corresponding interpolation figures.}
  \label{fig:RANK2 material}
\end{figure}

In this study, we consider the two-dimensional composite material, whcih is cross-shaped microstructure composed of two orthogonal rectangular bars intersecting at their centroid and rotated by an in-plane angle, as illustrated in Figure~\ref{fig:RANK2 material}(a). To globally parameterize the microstructures, we define a geometric variable set to

\begin{equation}
\mathbf{X} = \begin{bmatrix}
\mu_1^{(1)} & \mu_2^{(1)} & \theta^{(1)} \\
\mu_1^{(2)} & \mu_2^{(2)} & \theta^{(2)} \\
\vdots & \vdots & \vdots \\
\mu_1^{(N_\text{m})} & \mu_2^{(N_\text{m})} & \theta^{(N_\text{m})}
\end{bmatrix}
\in \mathbb{R}^{N_\text{m} \times m},
\end{equation}

where \(m\) is the number of control variable types (in here \(m=3\)), and each row \( \mathbf{X}^{(i)} = [\mu_1^{(i)}, \mu_2^{(i)}, \theta^{(i)}] \) represents the geometric parameters of the periodic microstructure assigned to the \( i \)-th macro element, with \( \mu_1 \) and \( \mu_2 \) denoting the widths of the horizontal and vertical bars, and \( \theta \) indicating their in-plane rotation angle.

Each unique parameter triplet \( \mathbf{X}^{(i)} \) defines a distinct microstructure with its corresponding homogenized elasticity tensor \( \mathbf{C}^{\text{H}}(\mu_1^{(i)}, \mu_2^{(i)}) \in \mathbb{R}^{3 \times 3} \), expressed in Voigt notation. To account for arbitrary orientation, the rotated effective tensor is computed as

\begin{equation}
\tilde{\mathbf{C}}^{\text{H}}(\mu_1^{(i)}, \mu_2^{(i)}, \theta^{(i)}) = \mathbf{T}^\top(\theta^{(i)}) \, \mathbf{C}^{\text{H}}(\mu_1^{(i)}, \mu_2^{(i)}) \, \mathbf{T}(\theta^{(i)}),
\end{equation}

where the transformation matrix \( \mathbf{T}(\theta) \in \mathbb{R}^{3 \times 3} \) is given by

\begin{equation}
\mathbf{T}(\theta) =
\begin{bmatrix}
\cos^2 \theta & \sin^2 \theta & 2 \sin \theta \cos \theta \\
\sin^2 \theta & \cos^2 \theta & -2 \sin \theta \cos \theta \\
-\sin \theta \cos \theta & \sin \theta \cos \theta & \cos^2 \theta - \sin^2 \theta
\end{bmatrix}.
\end{equation}

This matrix rotates the elasticity tensor from the local composite direction to the global simulation frame. Since the considered composite material exhibits orthotropic behavior, its homogenized elasticity tensor contains only four independent components: \( C_{11}^{\text{H}}, C_{22}^{\text{H}}, C_{12}^{\text{H}} \), and \( C_{33}^{\text{H}} \). These components are evaluated over a uniform grid in the geometric parameter space \( (\mu_1, \mu_2) \in [0.0, 1.0]^2 \). To enable efficient and continuous evaluation of material properties during the optimization process, each component \( C_{ij}^{\text{H}}(\mu_1, \mu_2) \) is approximated using a bivariate polynomial expansion of the form:

\begin{equation}
C^{\text{H}}_{ij}(\mu_1, \mu_2) \approx 
\sum_{q=0}^{q_{\max}} \sum_{r=0}^{r_{\max}} \alpha_{ij}^{(q,r)} (2\mu_1 - 1)^q (2\mu_2 - 1)^r,
\end{equation}

where \( \alpha_{ij}^{(q,r)} \) are the polynomial coefficients determined via least-squares fitting based on the precomputed homogenization dataset. This surrogate model facilitates rapid property queries for any macro element configuration \( \mathbf{X}^{(i)}\). The corresponding interpolation figures are shown in Figure~\ref{fig:RANK2 material}(b).

Based on the above content, we can explicitly establish a mapping between the geometric variables \( \mathbf{X} \) and the corresponding physical properties of the microstructures. To further reduce the total number of design variables, PCA is applied globally to compress the geometric parameter field \( \mathbf{X} \in \mathbb{R}^{N_\text{m} \times m} \). Specifically, the entire spatial distribution of geometric parameters is projected onto a low-dimensional subspace spanned by a small number of principal modes, thereby retaining the dominant geometric variations while significantly reducing computational complexity. Specifically, the reduced representation \( \mathbf{z} \in \mathbb{R}^{N_\text{c} \times m} \) is obtained as

\begin{equation}
\mathbf{z} = \mathbf{\Phi}^\top (\mathbf{X} - \bar{\mathbf{X}}),
\end{equation}

where \( \mathbf{\Phi} \in \mathbb{R}^{N_\text{m} \times N_\text{c}} \) contains the leading \( N_\text{c} \) principal components, and \( \bar{\mathbf{X}} \in \mathbb{R}^{N_\text{m} \times m} \) denotes the mean geometric field. The matrix \( \mathbf{z} \) captures the dominant variation modes and serves as the primary design variable during optimization. After optimization, the original geometric parameter field \( \mathbf{X} \) is approximately reconstructed through the inverse PCA mapping:

\begin{equation}
\mathbf{X} \approx \bar{\mathbf{X}} + \mathbf{\Phi} \mathbf{z}.
\end{equation}

This reduced-order modeling approach significantly decreases the number of design degrees of freedom from \( 3N_\text{m} \) to \( 3N_\text{c} \), while preserving the essential characteristics of the microstructure field. The choice of \( N_\text{c} \) controls the trade-off between model fidelity and computational efficiency.

Finally, the explicit geometry is obtained by \(\mathbf{X}\) through a geometry mapping operator:

\begin{equation}
\texttt{geom} = G_{\text{map}}(\mathbf{X}),
\end{equation}
  
where \texttt{geom} contains the explicit geometric representation, generally it is defined as a general CAD format file (like .dxf file), and will be used to achieve high-fidelity evaluation for optimization. The detail explanation of operator \( G_{\text{map}}(\cdot) \) will be introduced in the next subsection.

\begin{figure}[!tbp]
  \centering
  \includegraphics[width=\textwidth]{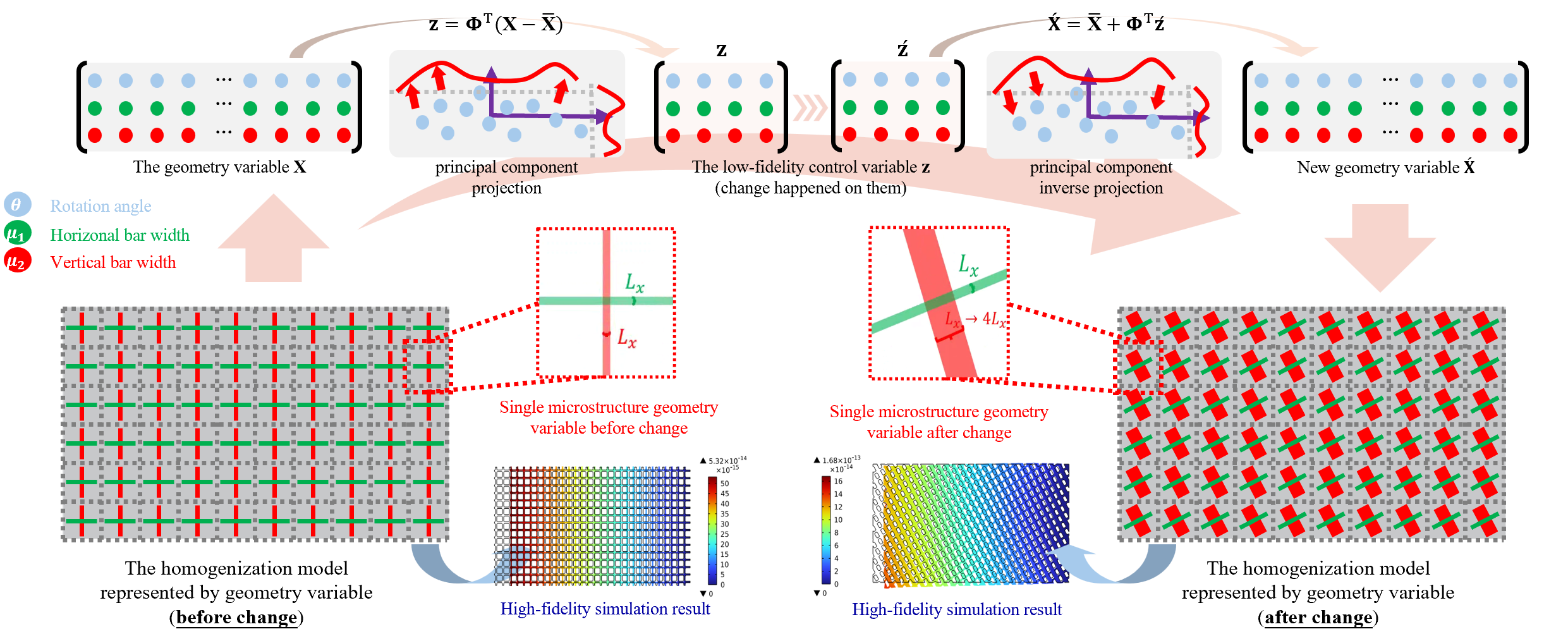}
  \caption{{The overview of the multifidelity modeling process.}}
  \label{fig:The cpt of hysp struct}
\end{figure}

This two-stage mapping—from the low-dimensional vector \( \mathbf{z} \) to the full-resolution geometry \texttt{geom}—facilitates efficient integration of compact design representations into high-fidelity simulations or CAD-based modeling workflows. The illustration of the workflow is shown in Figure~\ref{fig:The cpt of hysp struct}.

\subsection{Geometry mapping} \label{subsec:2.2}
A critical step in constructing a multifidelity model is geometry mapping, which transforms a finite set of low-dimensional geometric variables into explicit macrostructure geometry. In this work, geometry mapping consists of two major components: (1) converting the geometric variable \( \mathbf{X} \) into a detailed binary structure \(\boldsymbol{\phi}\) via de-homogenization method, and (2) extracting CAD-compatible model (\texttt{geom}) from the resulting structure field \(\boldsymbol{\phi}\) for downstream simulation or fabrication. The detail explaination for these components will be given in the following subsections.

\subsubsection{Detail pixal-based structure from phasor-based de-homogenization} \label{subsec:2.2.1}

The homogenization-based optimized design remains non-manufacturable due to scale disparity, necessitating a de-homogenization procedure that transforms theoretical multiscale configurations into physically realizable single-scale geometries. In this study, we adopt a phasor-based de-homogenization approach that transforms the geometric parameter field \( \mathbf{X} \in \mathbb{R}^{N_\text{m} \times 3} \) into a binary structural distribution \( \boldsymbol{\phi}(\mathbf{x}) \in \{0, 1\}^{N_\text{f}} \), where \( \mathbf{x} \in \mathbb{R}^2 \) represents spatial coordinates.

\begin{figure}[!tbp]
  \centering
  \includegraphics[width=\textwidth]{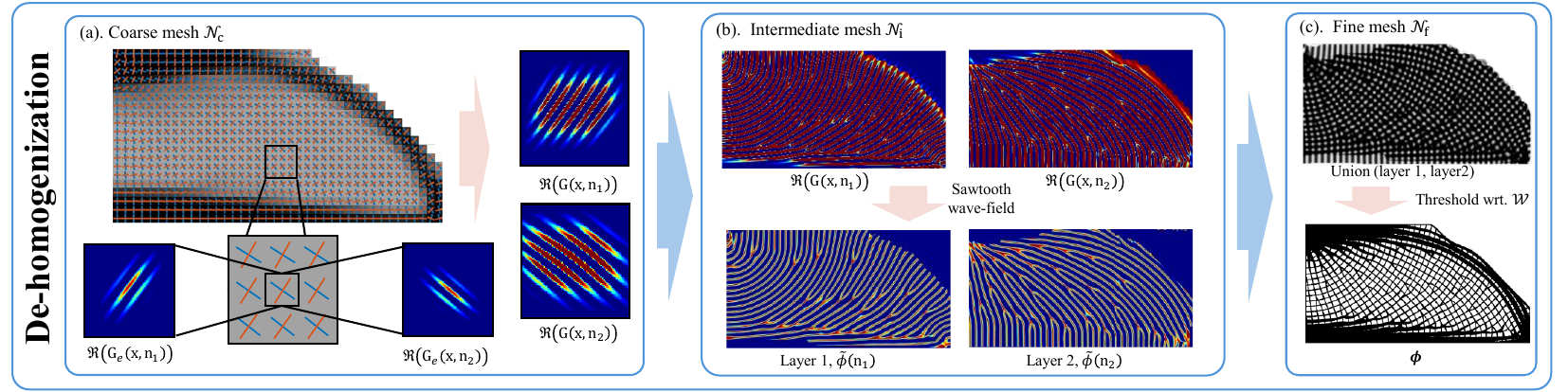}
  \caption{The mapping process for the phasor-based de-homogenization.}
  \label{fig:phasor_pipeline}
\end{figure}

As illustrated in Figure~\ref{fig:phasor_pipeline}, the de-homogenization process operates across three progressively refined computational domains: the coarse mesh \(\mathcal{N}_{c} \), intermediate mesh \(\mathcal{N}_{i} \), and fine mesh \( \mathcal{N}_{f} \).

On \(\mathcal{N}_{c} \), the parameter set \( \mathbf{X} \) contains the width parameters \( \mu_1, \mu_2 \) and the orientation angle \( \theta \) that characterize the orthotropic microstructure assigned to each macro element. From these parameters, two orthogonal lamination directions \( \mathcal{F} = \{ \mathbf{n}_1(\theta), \mathbf{n}_2(\theta) \} \) are defined, where \( \mathbf{n}_j(\theta) \) denotes the unit vector aligned with the \( j \)-th lamination direction. Correspondingly, directional thickness fields \( \mathcal{W} = \{ \overline{t}_{w1}(\mu_1), \overline{t}_{w2}(\mu_2) \} \) are introduced to control the layer widths.

In the intermediate domain \(\mathcal{N}_{i} \), oscillatory wave-fields are synthesized using phasor kernel operators. For each coarse element \( e \), a complex-valued phasor signal \( G_e(\mathbf{x}) \) is emitted from the element center \( \mathbf{x}_e \in \mathcal{N}_{c} \) along the local lamination direction \( \mathbf{n}_e \), and is sampled at the spatial location \( \mathbf{x} \in \mathcal{N}_{i} \) as

\begin{equation}
G_e(\mathbf{x}) = \exp\left( -\tilde{\beta} \Delta_e^{\mathbf{x}}(\mathbf{x}) \right) \exp\left( 2\pi i \omega \mathbf{n}_e \cdot (\mathbf{x} - \mathbf{x}_e) + i \varphi_e \right).
\end{equation}

Here, \( \tilde{\beta} \) is the spatial bandwidth parameter controlling the Gaussian decay perpendicular to the lamination direction; \( \Delta_e^{\mathbf{x}}(\mathbf{x}) \) is an anisotropic distance function defining the elliptical attenuation shape; \( \omega \) is the oscillation frequency governing the phasor wavelength; and \( \varphi_e \) is the phase offset assigned to each element \( e \)~\cite{woldseth2024phasor}.

The global wave-field \( G(\mathbf{x}, \mathbf{n}_j) \) corresponding to each principal direction \( \mathbf{n}_j \) is obtained by coherently superposing the filtered phasor emissions:

\begin{equation}
G(\mathbf{x}, \mathbf{n}_j) = \sum_{e} A_e(\mathbf{x}) G_e(\mathbf{x}, \mathbf{n}_j), \quad \mathbf{x} \in \mathcal{N}_{i},
\end{equation}

where \( A_e(\mathbf{x}) \) is a spatial sampling filter that smoothly regularizes the influence of each kernel. The real parts \( \Re(G(\mathbf{x}, \mathbf{n}_j)) \) are visualized in the intermediate mesh.

The directional laminate structures are extracted by applying a \texttt{sawtooth} wave demodulation to the synthesized fields, yielding the phase fields:

\begin{equation}
\tilde{\phi}(\mathbf{x}, \mathbf{n}_j) = \texttt{Sawtooth}\left( \text{Arg}(G(\mathbf{x}, \mathbf{n}_j)) \right), \quad j = 1,2,
\end{equation}

where \( \tilde{\phi}(\mathbf{x}, \mathbf{n}_j) \) represents the demodulated wave pattern aligned with the direction \( \mathbf{n}_j \).

In the final stage, the phase fields \( \tilde{\phi}(\mathbf{x}, \mathbf{n}_1) \) and \( \tilde{\phi}(\mathbf{x}, \mathbf{n}_2) \) are interpolated onto the fine mesh \( N_\text{f} \), thresholded according to the corresponding directional thickness fields \( \mathcal{W} \), and merged through Boolean union to generate the binary microstructure:

\begin{equation}
\phi(\mathbf{x}) = \mathcal{H}\left( \tilde{\phi}(\mathbf{x}, \mathbf{n}_1), \overline{t}_{w1} \right) \cup \mathcal{H}\left( \tilde{\phi}(\mathbf{x}, \mathbf{n}_2), \overline{t}_{w2} \right), \quad \mathbf{x} \in \mathcal{N}_{f},
\end{equation}

where \( \mathcal{H}(\cdot, \cdot) \) denotes a thresholding operator that discretizes continuous wave patterns into binary layers.

\subsubsection{Explicit geometry extraction and postprocessing} \label{subsec:2.2.2}

The binary field \( \boldsymbol{\phi}(\mathbf{x}) \) is pixel-based and thus unsuitable for CAD expression. To obtain a clean boundary representation, we extract the solid-void interface as an isocontour of \( \boldsymbol{\phi} \) at a prescribed threshold \( c \), yielding a polyline boundary \( \mathcal{P} \subset \mathbb{R}^2 \). This process corresponds to the internal geometry representation \( \texttt{geom} \), which includes nodal coordinates and topological information.

\subsubsection*{Marching squares algorithm}

The isocontour corresponding to \( \boldsymbol{\phi}(\mathbf{x}) = c \) is extracted using the Marching Squares Algorithm (MSA)~\cite{lorensen1998marching}, where \( c \) denotes the prescribed threshold for boundary extraction (typically set to 0.5 for binary structures).

Each grid cell in the structured mesh is defined by the scalar field values at its four corners, denoted as \( \phi_{00}, \phi_{10}, \phi_{01}, \phi_{11} \), corresponding respectively to the bottom-left, bottom-right, top-left, and top-right nodes. These scalar values are first binarized according to a prescribed threshold \(c\):

\begin{equation}
b_k =
\begin{cases}
1, & \phi_k \geq c, \\[4pt]
0, & \phi_k < c,
\end{cases}
\quad k \in \{00, 10, 01, 11\}.
\end{equation}

Each resulting binary configuration forms a 4-bit code that maps the cell into one of 16 predefined intersection cases, which determine how the isocontour intersects the edges of the cell. For each edge where the binary states of its two end nodes differ (i.e., one node is inside and the other outside the thresholded region), the precise interface point \( \mathbf{p}_{\text{int}} \) is computed by linear interpolation between the corresponding corner points \( \mathbf{p}_{\text{end},1} \) and \( \mathbf{p}_{\text{end},2} \), whose scalar values are \( \phi_1 \) and \( \phi_2 \), respectively:

\begin{equation}
\mathbf{p}_{\text{int}} = \mathbf{p}_{\text{end},1} + \frac{c - \phi_1}{\phi_2 - \phi_1} (\mathbf{p}_{\text{end},2} - \mathbf{p}_{\text{end},1}).
\end{equation}

All such interpolated points are collected to form the interface polyline:

\begin{equation}
\mathcal{P} = \left\{ \mathbf{p}_{\text{int}}^{(1)}, \mathbf{p}_{\text{int}}^{(2)}, \dots, \mathbf{p}_{\text{int}}^{(N_\text{msa})} \right\}, \quad \mathbf{p}_{\text{int}}^{(i)} \in \mathbb{R}^2,
\end{equation}

where \(N_\text{msa}\) is the number of control points after MSA operation, which accurately traces the interface between the solid and void regions.

\subsubsection*{Equal arc-length resampling}

To enable robust CAD modeling and meshing, the non-uniformly distributed polyline \( \mathcal{P} = \{ \mathbf{p}_{\text{int}}^{(i)} \} \) is resampled to obtain a uniformly spaced point set. We use a method named Equal Arc-length Resampling (EAR) to achieve this.

Firstly, the cumulative arc length along \( \mathcal{P} \) is computed as

\begin{equation}
s_i = \sum_{j=1}^{i-1} \|\mathbf{p}_{\text{int}}^{(j+1)} - \mathbf{p}_{\text{int}}^{(j)}\|, \quad s_1 = 0, \quad s_{N_\text{msa}} = L,
\end{equation}

where \( L \) denotes the total length of the original polyline. The target number of resampling points \(N_\text{ear}\) is control by a simplication factor \( M \), where \(N_\text{ear} = M \cdot\ N_\text{msa}\), the uniform arc-length interval is defined as

\begin{equation}
\Delta s = \frac{L}{N_\text{ear}-1}.
\end{equation}

For each resampled point \( k \in \{1, \dots, N_\text{ear}\} \), the new coordinate \( \tilde{\mathbf{p}}_k \) is obtained by linearly interpolating between two adjacent original points \( \mathbf{p}_{\text{int}}^{(i)} \) and \( \mathbf{p}_{\text{int}}^{(i+1)} \) such that \( s_i \leq (k-1)\Delta s < s_{i+1} \):

\begin{equation}
\tilde{\mathbf{p}}_k = \mathbf{p}_{\text{int}}^{(i)} + \frac{(k-1)\Delta s - s_i}{s_{i+1} - s_i}  \left( \mathbf{p}_{\text{int}}^{(i+1)} - \mathbf{p}_{\text{int}}^{(i)} \right).
\end{equation}

Collecting all resampled points yields the uniformly distributed polyline:

\begin{equation}
\tilde{\mathcal{P}} = \left\{ \tilde{\mathbf{p}}_1, \tilde{\mathbf{p}}_2, \dots, \tilde{\mathbf{p}}_{N_\text{ear}} \right\}.
\end{equation}

\subsubsection*{geometric deviation evaluation}

To evaluate the shape fidelity between the original and resampled polylines, the geometric deviation is defined as
\begin{equation}
\varepsilon_{\text{shape}} = \max \left\{
    \max_{\tilde{\mathbf{p}} \in \tilde{\mathcal{P}}} \min_{\mathbf{p}_{\text{int}} \in \mathcal{P}} \|\tilde{\mathbf{p}} - \mathbf{p}_{\text{int}}\|,
    \max_{\mathbf{p}_{\text{int}} \in \mathcal{P}} \min_{\tilde{\mathbf{p}} \in \tilde{\mathcal{P}}} \|\mathbf{p}_{\text{int}} - \tilde{\mathbf{p}}\|
\right\},
\end{equation}

where the first term measures the maximum distance from resampled points to the original curve, and the second term measures the reverse. This deviation is then normalized by the original arc length \( L \) to yield the relative geometric accuracy:

\begin{equation}
\varepsilon_{\text{acc}} = \left(1 - \frac{\varepsilon_{\text{shape}}}{L} \right) \times 100\%.
\end{equation}

\subsubsection*{CAD export and utilization}

The final resampled polyline \( \tilde{\mathcal{P}} \) could be directly exported into a \texttt{.dxf} format, encoding the point coordinates and connectivity in a CAD-compatible structure. These files are subsequently imported into commercial CAE software, where the closed contours are meshed for high-fidelity finite element simulations.

During the meshing phase, a third-party mesh generation software (such as COMSOL Multiphysics~\cite{multiphysics1998introduction}) is employed to discretize the extracted geometry using free triangular elements. To ensure controllability and consistency, the minimum edge length of the generated triangular mesh is set equal to the arc-length interval \( \Delta s \) defined during the resampling process. As a result, by adjusting the control parameter \( M \), which specifies the number of resampled boundary points, both the geometric accuracy and the final number of mesh elements can be directly regulated.

\begin{figure}[!tbp]
  \centering
  \includegraphics[width=\textwidth]{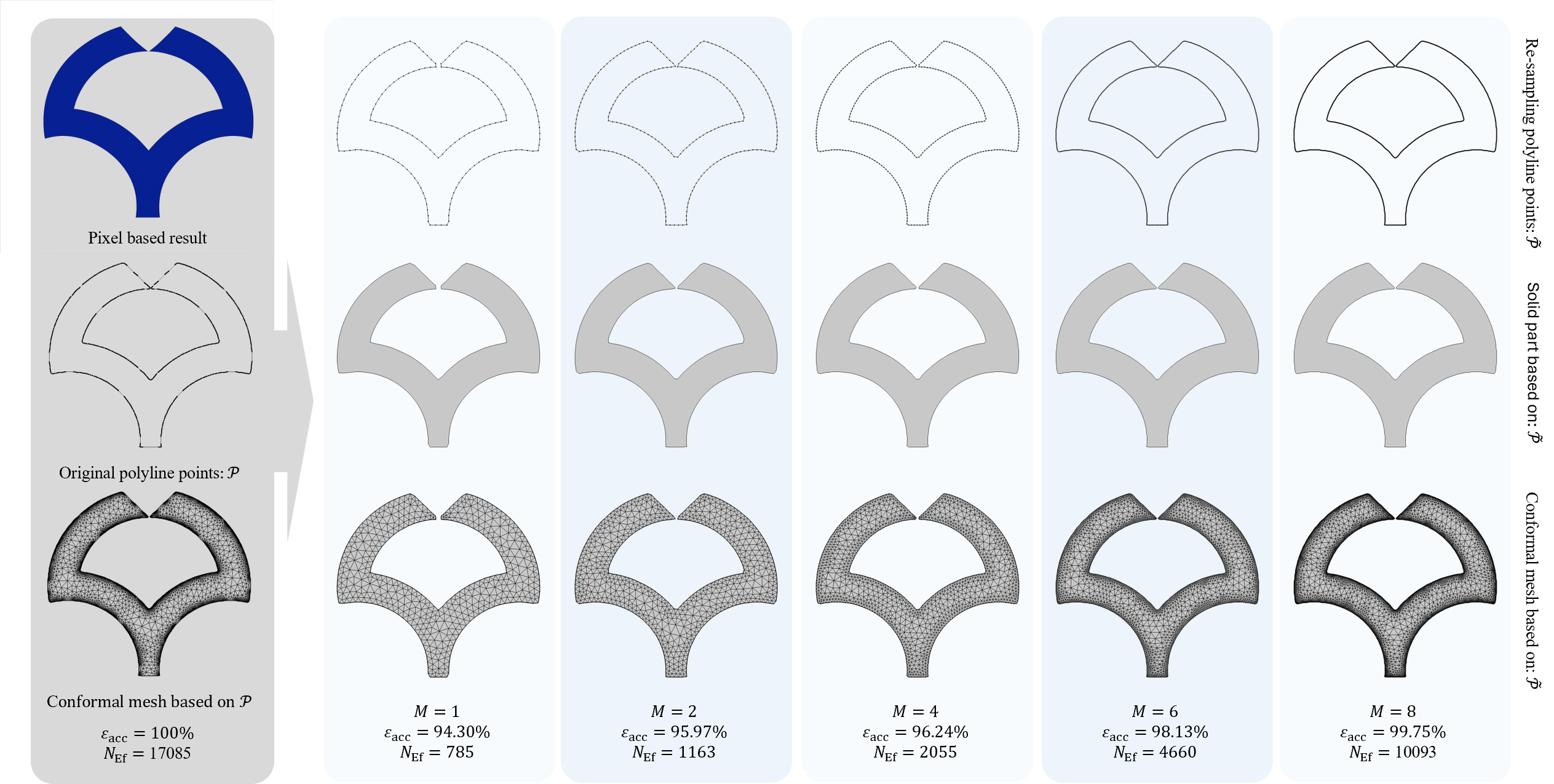}
  \caption{{Schematic illustration of geometry extraction and postprocessing.}}
  \label{fig:The CAD}
\end{figure}

As illustrated in Figure~\ref{fig:The CAD}, we demonstrate this process using a given binary image. The top row shows the mesh generated directly from the MSA without resampling. Due to the non-uniform distribution of boundary points, the resulting mesh exhibits significant irregularity, with densely packed elements in some regions and sparse elements in others. In contrast, the bottom three rows present meshes generated using the EAR method for different values of \( M \). It can be observed that, as \( M \) increases, the geometric fidelity of the mesh improves, the boundary approximation becomes more accurate, and the element distribution becomes more uniform. Additionally, compared to the direct MSA-based extraction, the EAR-based meshes achieve a higher quality discretization with fewer overall elements.

In concluison, the transformation pipeline \( G_{\text{map}}: \mathbf{X} \rightarrow \boldsymbol{\phi}(\mathbf{x}) \rightarrow \mathcal{P} \rightarrow \tilde{\mathcal{P}} \rightarrow \texttt{geom} \) establishes a complete and manufacturable link between abstract geometric control variables and explicit structural realizations.

\section{Optimization model and its numerical implementation} \label{sec:Optimization model}
\subsection{Optimization formulation with latent space evolution}
As discussed previously, the transformation from the reduced control variable \( \mathbf{z} \) to full-resolution geometry \texttt{geom} is inherently irreversible, as the mapping operator \( G_{\text{map}}(\cdot) \) involves heuristic and non-differentiable operations, whcih hinder the application of gradient-based optimization methods due to the lack of explicit sensitivity information. 
Furthermore, the discrepancy in information fidelity between control variables and high-fidelity simulations introduces strong nonlinearity to the design landscape. To address these challenges, we employ a data-driven approach based EA. In EA, maintaining population diversity is critical to avoid premature convergence and to explore a broad design space. We quantify design diversity through the \textit{volume fraction} (i.e., material usage) of the high-fidelity geometry, denoted by \( G_{\mathrm{vf}} \). However, structural performance \( G_{\mathrm{opt}} \) typically improves with increased material usage, introducing a natural trade-off. To maintain this balance, we adopt a multi-objective optimization strategy that considers both metrics simultaneously. The optimization problem is formulated as

\begin{equation}
\begin{array}{ll}
\text{Find:} & \mathcal{J} = \{ \mathbf{z}_1, \mathbf{z}_2, \ldots, \mathbf{z}_n \} \\
\text{Minimize:} & \mathbf{O}(\mathbf{z}_i) = \left[ G_{\mathrm{vf}}(\mathbf{z}_i),\; G_{\mathrm{opt}}(\mathbf{z}_i) \right], \quad i = 1,\ldots, n \\
\text{Subject to:} & \mathbf{X}_i \in \Omega_c
\end{array}
\label{eq:latent-optimization}
\end{equation}

\begin{figure}[!tbp]
  \centering
  \includegraphics[width=\textwidth]{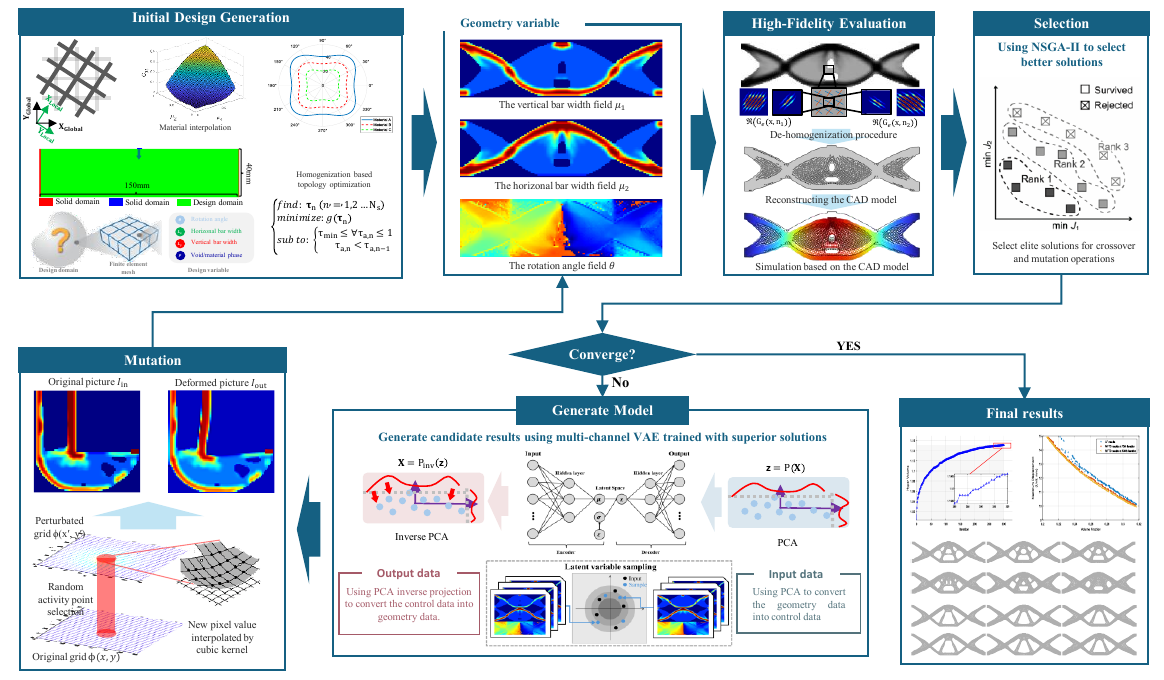}
  \caption{{Schematic illustration of geometry extraction and postprocessing.}}
  \label{fig:The MFTD flow}
\end{figure}

Here, \(\mathcal{J}\) represents the design population, \(n\) is the number of sampling, \(\mathbf{z}_i \in \mathbb{R}^{N_c \times 3}\) is the reduced dimension design vector of the \( i \)-th individual obtained from the PCA operation, and \(\Omega_c\) is the feasible design space defined by physical and geometrical constraints. The overall optimization framework consists of five key components (Figure~\ref{fig:The MFTD flow})

\begin{itemize}
    \item \textbf{Initial Design Generation:} Low-fidelity initial designs are generated via homogenization-based topology optimization (which could be also called as low-fidelity optimization), providing physically feasible initial solutions that satisfy fundamental constraints.
    
    \item \textbf{High-Fidelity Evaluation:} Each latent design is mapped to a detailed geometry, followed by finite element simulation to obtain \(G_{\mathrm{vf}}\) and \(G_{\mathrm{opt}}\).
    
    \item \textbf{Selection:} The NSGA-II with elitism is used to retain promising candidates and maintain diversity in the population.
    
    \item \textbf{Crossover:} New designs are generated in the latent space through crossover, guided by a VAE-based generative model and specific latent space resampling method.
  
    \item \textbf{Mutation:} The mutation operation is implemented through an image processing technique, which efficiently introduces novel geometric variations and facilitates broader design space exploration with low computational cost.
\end{itemize}

\subsection{Initial design generation}
\label{subsec:3.1}
In non-gradient-based optimization algorithms, the quality of the initial population plays a crucial role in influencing both convergence speed and solution quality. To generate effective initial designs, we adopt an asymptotic homogenization-based topology optimization strategy. The general low-fidelity optimization problem is formulated as

\begin{equation}
\begin{aligned}
\text{Find:} \quad & \mathbf{X} \\
\text{Minimize:} \quad & \text{obj}(\mathbf{X}) \\
\text{Subject to:} \quad & \mathbf{K}(\mathbf{X}) \mathbf{U} = \mathbf{F} \\
& V(\mathbf{X}) \leq V_0 \\
& \mathbf{X}_{\min} \leq \mathbf{X} \leq \mathbf{X}_{\max}
\end{aligned}
\label{eq:low-fidelity}
\end{equation}

where \(\text{obj}(\cdot)\) denotes the objective function. The constraints enforce structural equilibrium, a global volume fraction limit, and geometric bounds. In this work, the low-fidelity optimization model only considers stiffness optimization under linear elasticity assumptions. \(\mathbf{K}\) is global stiffness matrix, \(\mathbf{U}\) is the corresponding nodal response field, and \(\mathbf{F}\) is the external force. The material peroperties of the microstructure in macroscale follow the strategy introduced in Subsection~\ref{subsec:2.1}, where the base material is assumed to be a pseudo-isotropic material with a Young's modulus of 1 and a Poisson's ratio of 0.3.

Depending on the specific design problem, the detailed construction of \(\mathbf{X}\), \(\mathbf{X}_{\min}\), \(\mathbf{X}_{\max}\) and the corresponding material interpolation schemes are adapted accordingly. The specific formulations for each case study are provided in Section~\ref{sec:Numerical Examples}.

\subsection{Selection}
\label{subsec:3.2}

Selection is responsible for screening and retaining high-quality designs during each evolutionary iteration. We adopt the elitism-based selection mechanism of NSGA-II, a widely used multi-objective evolutionary algorithm. This method ranks individuals based on Pareto dominance and ensures both convergence and diversity by preserving elite solutions and encouraging exploration.

\subsection{Crossover}
\label{subsec:3.3}

Crossover is a key evolutionary operation that combines genetic material from multiple parents to generate new candidate solutions. In the context of MFTD, this operation is performed in the latent space using a VAE. Since the VAE must be retrained at each optimization iteration based on the current elite samples, a compact and efficient VAE framework is critical to the successful implementation of MFTD.

In our previous work on multi-variable optimization~\cite{kawabe2025data, xu2025evolutionary}, each low-fidelity design sample was represented as a matrix \(\mathbf{X} \in \mathbb{R}^{N_m \times m}\), where \(N_m\) is the data dimension of each sampling point (analogous to pixels number) and \(m\) denotes the number of design field per point (analogous to color channels). Due to the absence of dimensionality reduction, convolutional VAEs was employed to process such high-dimensional multi-channel data, resulting in high computational costs during retraining.

Since we introduce a dimensionality reduction step using PCA in this work. For each sample, PCA will apply along the spatial dimension to reduce the number of sampling points from \(N_m\) to \(N_c\) while preserving the original \(m\) channels. The resulting reduced design variable field is represented as \(\mathbf{z} \in \mathbb{R}^{N_c \times m}\), where \(N_c \ll N_m\). Subsequently, \(\mathbf{z}\) is flattened into a one-dimensional vector of size \((m \times N_c) \times 1\), enabling the use of a simplified fully connected VAE architecture with two hidden layers in both the encoder and decoder. Although a simple fully connected VAE architecture is employed after flattening, the relatively low dimensionality of each individual sample ensures that the VAE retains sufficient capacity to capture inter-channel correlations. As a result, the trained VAE can successfully generate latent representations that maintain the intrinsic relationships between different channels, producing realistic and compliant design samples. For the architecture of the VAE, the encoder outputs a latent distribution defined by mean and variance vectors, \(\boldsymbol{\mu}, \boldsymbol{\sigma} \in \mathbb{R}^{N_{\text{lat}}}\) (\({N_{\text{lat}}}\) is the dimension of VAE latent space), from which a latent code \(\mathbf{z}_{\text{lat}}\) is sampled as

\begin{equation}
\mathbf{z}_{\text{lat}} = \boldsymbol{\mu} + \boldsymbol{\sigma} \circ \boldsymbol{\epsilon}, \quad \boldsymbol{\epsilon} \sim \mathcal{N}(\mathbf{0}, \mathbf{I}),
\label{eq:vae-sampling}
\end{equation}

where \(\boldsymbol{\epsilon} \sim \mathcal{N}(\mathbf{0}, \mathbf{I})\) denotes a standard Gaussian noise vector, introduced via the reparameterization trick to enable backpropagation through the stochastic sampling process. Then, the VAE is trained by minimizing the following loss function

\begin{equation}
\mathcal{L} = \mathcal{L}_{\text{recon}} + \beta\mathcal{L}_{\text{KL}},
\label{eq:vae-loss}
\end{equation}

where \(\mathcal{L}_{\text{recon}}\) is the reconstruction loss (typically mean squared error) and \(\mathcal{L}_{\text{KL}}\) is the Kullback–Leibler divergence

\begin{equation}
\mathcal{L}_{\text{KL}} = -\frac{1}{2} \sum_{i=1}^{N_{\text{lat}}} \left( 1 + \log \sigma_i^2 - \mu_i^2 - \sigma_i^2 \right).
\label{eq:kl-divergence}
\end{equation}

This probabilistic formulation promotes smooth interpolation and generalization. Once trained, the VAE enables efficient crossover operations by recombining latent representations within a compact and regularized space~\cite{kawabe2025data}.

\subsection{Mutation: deformation-based mutation via radial sampling perturbation}
\label{subsec:3.4}

Mutation introduces stochastic variations into the design population, promoting exploration and preventing premature convergence. In our framework, mutation is implemented via a differentiable image deformation process. Given an original image \( I_{\text{in},\mathrm{m}}(i, j) \in \mathbb{R} \), where \( (i, j) \in \mathbb{Z}^2 \) denote discrete pixel indices, the initial sampling coordinates \(\boldsymbol{\phi}_{\mathrm{m}}(i,j) = (x(i,j),\; y(i,j))\) are aligned with the regular grid and centered at integer locations. After deformation, we obtain a new set of sampling coordinates \(\boldsymbol{\phi}'_{\mathrm{m}}(i,j) = (x'(i,j),\; y'(i,j))\). Our goal is to construct a new image \( I_{\text{out},\mathrm{m}} \) by sampling from the original image \( I_{\text{in},\mathrm{m}} \), such that
\begin{equation}
I_{\text{out},\mathrm{m}}(\boldsymbol{\phi}_{\mathrm{m}}(i,j)) = I_{\text{in},\mathrm{m}}(\boldsymbol{\phi}'_{\mathrm{m}}(i,j))
\end{equation}

Since \( \boldsymbol{\phi}'_{\mathrm{m}}(i,j) \) may lie at non-integer locations, we use bicubic interpolation to approximate the value
\begin{equation}
I_{\text{out},\mathrm{m}}(\boldsymbol{\phi}_{\mathrm{m}}(i,j)) \approx 
\sum_{m=-1}^{2} \sum_{n=-1}^{2} 
I_{\text{in},\mathrm{m}}(\boldsymbol{\phi}_{\mathrm{m}}(i+m, j+n)) \omega\big(\boldsymbol{\phi}'_{\mathrm{m}}(i,j) - \boldsymbol{\phi}_{\mathrm{m}}(i+m, j+n)\big)
\end{equation}

where \( \omega(\cdot) \) is the Keys interpolation kernel function~\cite{keys2003cubic}, which controls the weight contribution of neighboring points to the final result, given by
\begin{equation}
w_c(r) =
\begin{cases}
( a + 2 )|r|^3 - ( a + 3 )|r|^2 + 1, & 0 \leq |r| < 1, \\
a|r|^3 - 5a|r|^2 + 8a|r| - 4a, & 1 \leq |r| < 2, \\
0, & |r| \geq 2,
\end{cases}
\label{eq:cubic_kernel}
\end{equation}

where \( a \in [-1, 0] \) is a smoothness parameter. In this work, we adopt \( a = -0.5 \). This interpolation ensures smooth transitions and second-order continuity, as each output value depends on a \( 4 \times 4 \) neighborhood in the original image grid.

\begin{figure}[!tbp]
  \centering
  \includegraphics[width=\textwidth]{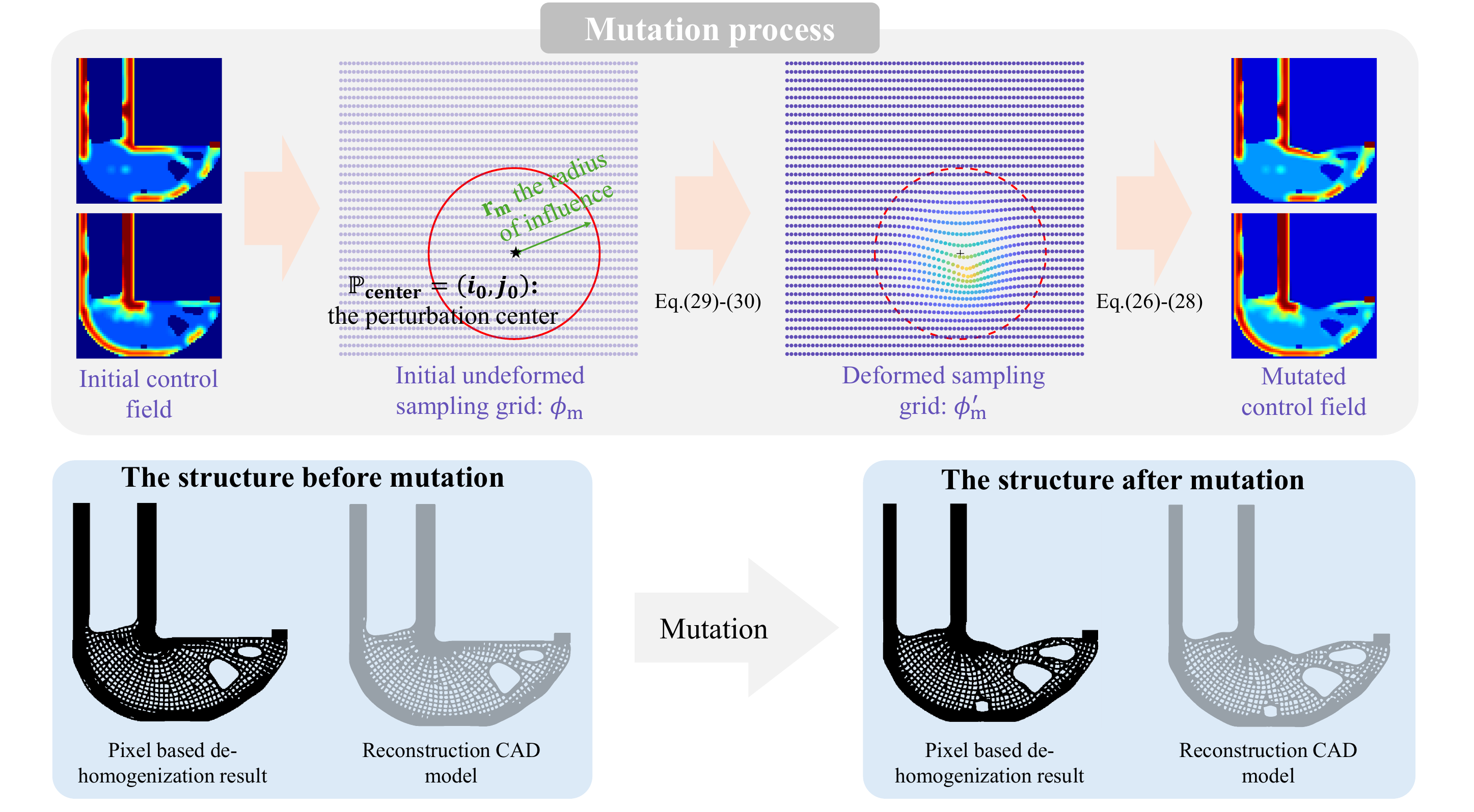}
  \caption{The whole process for the mutation method proposed in this work.}\label{fig:The mutation}
\end{figure}

Based on the above discussion, it is evident that the key mechanism behind mutation lies in the deformation of the interpolation grid. To this end, we define a radial deformation-based mutation operator, which perturbs the sampling coordinates according to a spatially decaying displacement field. Specifically, the deformed sampling coordinate \( \boldsymbol{\phi}'_{\mathrm{m}}(x, y) \) is defined as
\begin{equation}
\boldsymbol{\phi}'_{\mathrm{m}} = 
\left(
x - \alpha_{\mathrm{m}} u_{x,\mathrm{m}},\ 
y + \alpha_{\mathrm{m}} u_{y,\mathrm{m}}
\right),
\label{eq:phi_radial}
\end{equation}

where \( u_{x,\mathrm{m}} \in \mathbb{R}^{i \times j} \) and \( u_{y,\mathrm{m}} \in \mathbb{R}^{i \times j} \) are the components of a displacement field derived from a low-fidelity structural simulation, approximating the deformation response under a prescribed load. The attenuation factor \( \alpha_{\mathrm{m}} \in \mathbb{R}^{i \times j} \) modulates the influence of this field and is defined as
\begin{equation}
\alpha_{\mathrm{m}}(i, j) = 
\max\left(0,\ 1 - \frac{ \| (i, j) - (i_0, j_0) \| }{r_m} \right)^{p_0},
\label{eq:alpha_decay}
\end{equation}

where \( (i_0, j_0) \in \mathbb{R}^2 \) is the perturbation center, \( r_m \in \mathbb{R}^+ \) denotes the radius of influence, and \( p_0 \in \mathbb{R}^+ \) controls the smoothness of the decay. The resulting deformation enables a controllable and spatially localized mutation of the image structure, while preserving smoothness through interpolation. The whole mutation process is shown in Figure~\ref{fig:The mutation}.

\section{Numerical examples} \label{sec:Numerical Examples}
In this study, we used COMSOL Multiphysics (version 6.2), which is a commercial software using the Finite Element Method (FEM), for achieving the mesh discretion and solving the high-fidelity evaluation. Besides, VAE was implemented using Python (version 3.7.4) and TensorFlow (version 2.0.0). The other part of the algorithm given in section 3 was manipulated using MATLAB (version 2022b). All simulations and computations in this work were conducted on a high-performance workstation equipped with dual AMD EPYC 7763 CPUs (128 cores in total, 2.45-3.5 GHz), 1 TB of RAM, and two NVIDIA RTX A6000 GPUs.

\subsection{Double clamper beam design} \label{sec:4.1}
The first 2D example involves designing a double-clamped beam, serving as a fundamental benchmark for exploring various key issues in this case. The boundary conditions and structural dimensions are illustrated in Figure~\ref{fig:The bc of db}(a). The beam has a rectangular shape, with its width greater than its length. The degrees of freedom on the left and right sides are fully fixed, restricting movement in the horizontal and vertical directions. A horizontally downward force is applied at the top center. In this case, two design objectives will be considered: 1. Stiffness optimization, and 2. Buckling performance improvement under a stiffness constraint.

\subsubsection{Preparation of initial design} \label{sec:4.1.1}
The low-fidelity optimization model follows the general formulation in Eq.~\eqref{eq:low-fidelity}. But to further enhance material distribution at the macroscopic scale, we introduce an additional macroscopic material variable \( \rho \in [0,1] \) to distinguish between void regions and lattice-filled regions. Accordingly, the effective elasticity tensor is modified to incorporate the influence of \( \rho \) as
\begin{equation}
\tilde{\mathbf{C}}^{\text{H}}(\mu_1, \mu_2, \theta, \rho) = \mathbf{T}^\top(\theta) \, \mathbf{C}^{\text{H}}(\mu_1 \rho, \mu_2 \rho) \, \mathbf{T}(\theta),
\label{eq:macro-rho}
\end{equation}

The geometric variable vector is now defined as
\begin{equation}
  \mathbf{X} = \begin{bmatrix} \mu_1 \rho & \mu_2 \rho & \theta \end{bmatrix}^\top.
\end{equation}

For the low-fidelity optimization model, the specified optimization model is formulated as

\begin{equation}
  \begin{array}{ll}
  \text{Find:} & \mu_1,\, \mu_2,\, \theta, \rho \\
  \text{Minimize:} & \text{obj}(\mu_1, \mu_2, \theta, \rho) \\
  \text{Subject to:} & \mathbf{K}(\mu_1, \mu_2, \theta, \rho) \mathbf{U} = \mathbf{F} \\
  & V(\mu_1, \mu_2, \rho) \leq V_0 \\
  & L_{\min} \leq \mu_1,\, \mu_2 \leq L_{\max} \\
  & \rho_{\min} \leq \rho \leq \rho_{\max} \\
  & -4\pi \leq \theta \leq 4\pi
  \end{array}
\end{equation}

\begin{figure}[!tbp]
  \centering
  \includegraphics[width=\textwidth]{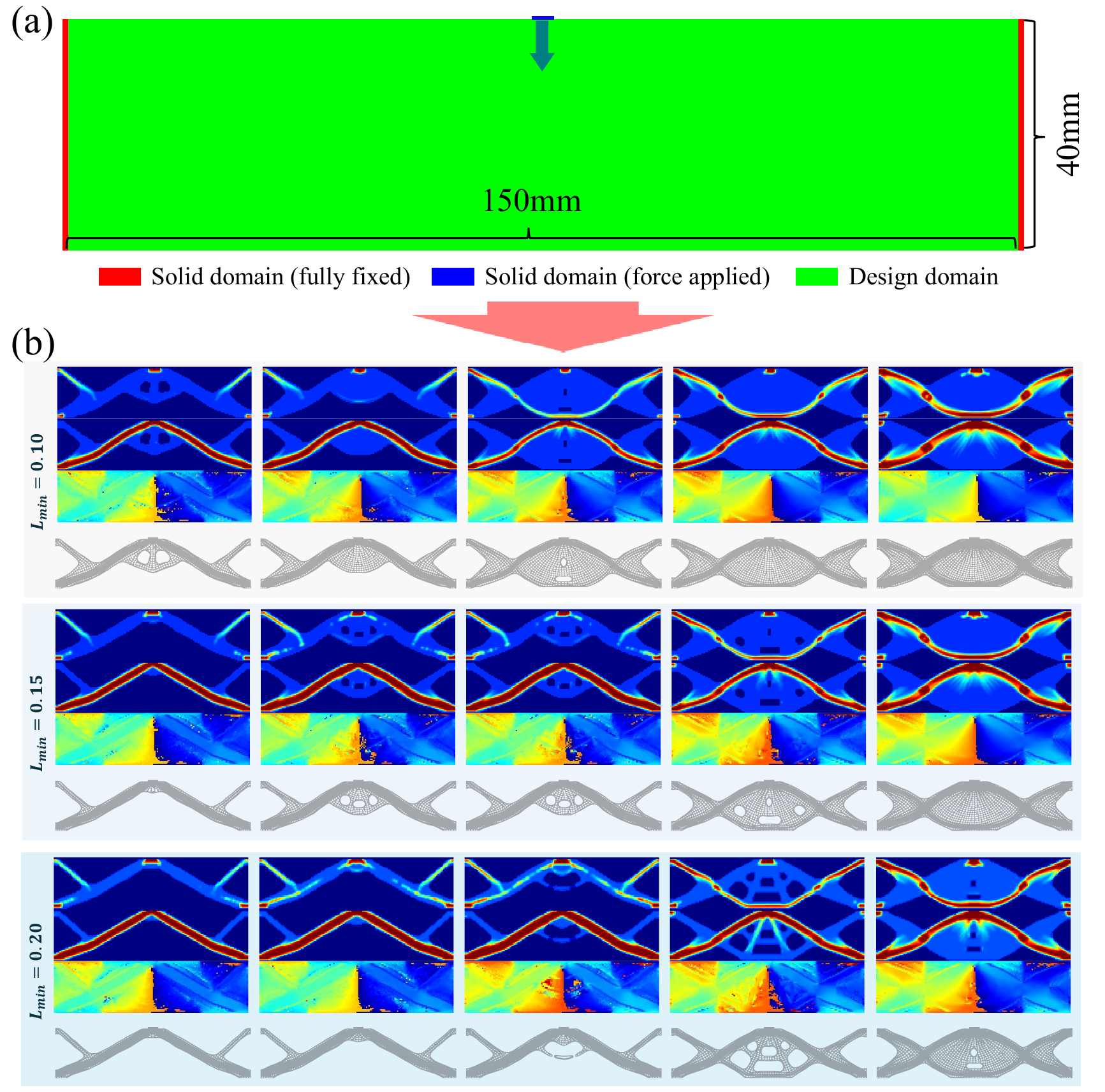}
  \caption{(a) The boundary condition and the design space of double clamper beam; (b) the 15 representative initial solutions with different \(L_{\min}\).}\label{fig:The bc of db}
\end{figure}

To ensure diversity in the initial solutions, we achieve this by controlling the volume fraction \(V_0\) and the corresponding minimum size \(L_{\min}\). Specifically, we generate 160 different solutions by varying the \(V_0\) from 0.25 to 0.50. Subsequently, for each case, we apply three different minimum feature sizes (\(L_{\min}\) = 0.1, 0.15, and 0.2), resulting in a total of 480 solutions. In Figure~\ref{fig:The bc of db}(b), we present 15 representative initial solutions, and their CAD model after geometry mapping. The selection includes five structures for each \(L_{\min}\), arranged in increasing volume fraction order. Observations indicate that this approach effectively ensures structural topology diversity in the initial solutions. For example, when \(L_{\min}\) is relatively large, more large internal voids tend to appear.

\subsubsection{Stiffness improvement} \label{sec:4.1.2}

\begin{figure}[!tbp]
  \centering
  \includegraphics[width=\textwidth]{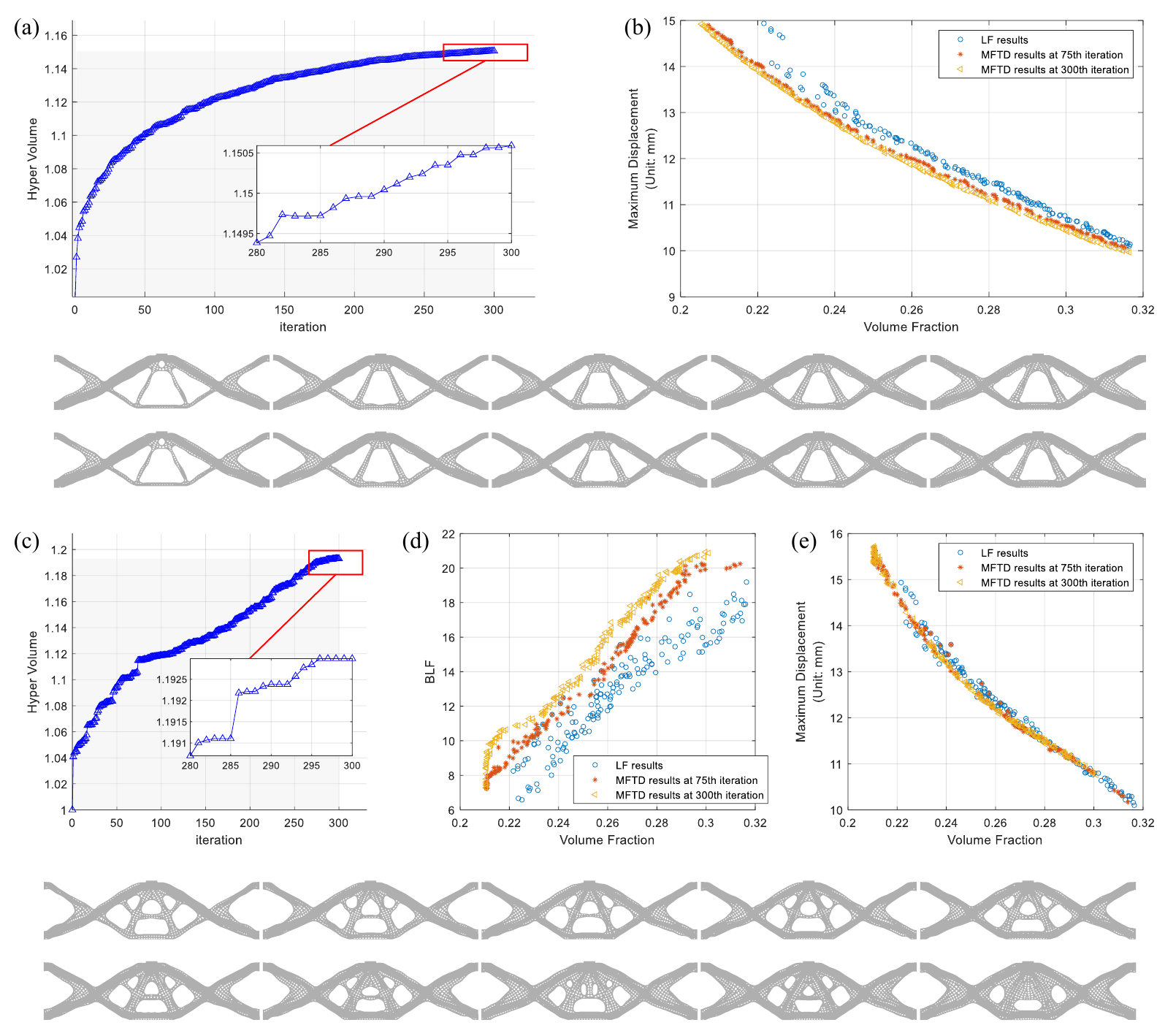}
  \caption{The optimization result for the double clamper beam, (a)--(b) are the results for the stiffness based design, and (c)--(e) are the results for the buckling based design. (a) and (c) are the hypervolume value convergence history, and (b)(d)(e) are the corresponding pareto front of the responses.}\label{fig:The result of db}
\end{figure}

Based on the initial solutions shown in Subsection 4.1.1, we first optimize the structure with the objective of minimizing maximum displacement and volume. The optimization process consists of 300 iterations, and the corresponding hypervolume convergence curve is shown in Figure~\ref{fig:The result of db}(a). Compared to the initial solution, the final optimized hypervolume shows a 15.1\% improvement. 

The Pareto front of Maximum Displacement vs. Volume is plotted in Figure~\ref{fig:The result of db}(b). As the volume increases, the maximum displacement gradually decreases, demonstrating a typical optimization trade-off. This trend indicates that increasing the structural volume reduces maximum displacement, but at the cost of higher material consumption. At the same volume, the MFTD method generally achieves lower maximum displacement (with the curve shifting downward) and provides solutions closer to the Pareto optimal front. This suggests that the MFTD method identifies superior design solutions, allowing the structure to maintain low maximum displacement with a smaller volume. 

Subsequently, we selected 10 representative MFTD-optimized structures based on their volume distribution. In terms of stiffness design, from a macroscopic material distribution perspective, the final result is highly similar to traditional SIMP, and large-area lattice filling is minimized as much as possible.

\subsubsection{Buckling improvement with stiffness constraint} \label{sec:4.1.3}

Linear buckling performance is considered in this case, starting with the same initial design as the stiffness-based case. Here, we select the Buckling Load Factor (BLF) of the first buckling mode and its corresponding volume as evaluation metrics. Additionally, a stiffness constraint is imposed, ensuring that the optimized structure's stiffness is not lower than that of the initial design with the same volume. To unify the evaluation, maximum displacement is used as a proxy for stiffness. Similarly, the optimization process consists of 300 iterations, and the corresponding hypervolume convergence curve is shown in Figure~\ref{fig:The result of db}(c). Compared to the initial solution, the final optimized hypervolume shows a 19.4\% improvement. The Pareto front of BLF vs. Volume is plotted in Figure~\ref{fig:The result of db}(d), and the Pareto front of Maximum Displacement vs. Volume is shown in Figure~\ref{fig:The result of db}(e). 

Similar to the observations in stiffness optimization, as volume increases, the maximum BLF also increases. At the same volume, the MFTD method consistently achieves higher BLF values (with the curve shifting upward) and provides solutions closer to the Pareto optimal front. For maximum displacement, however, the distribution of the MFTD method remains largely consistent with the initial solution, without noticeable upward or downward shifts. This indicates that the stiffness constraint has been effectively enforced, ensuring that the structure does not sacrifice overall stiffness performance and stability in pursuit of a higher BLF value in the first mode. 

Subsequently, we selected 10 representative MFTD-optimized structures based on their volume distribution. In terms of buckling design, more geometric features appear to enhance structural stability. However, since our initial solution is based on stiffness, and strong stiffness constraints are imposed during the MFTD process, the optimization does not result in a fully buckling-optimized structure. Instead, it selectively improves buckling performance while minimizing the impact on overall stiffness.

\subsubsection{Comprehensive analysis for the double clamper beam design} \label{sec:4.1.4}

\begin{figure}[!tbp]
  \centering
  \includegraphics[width=\textwidth]{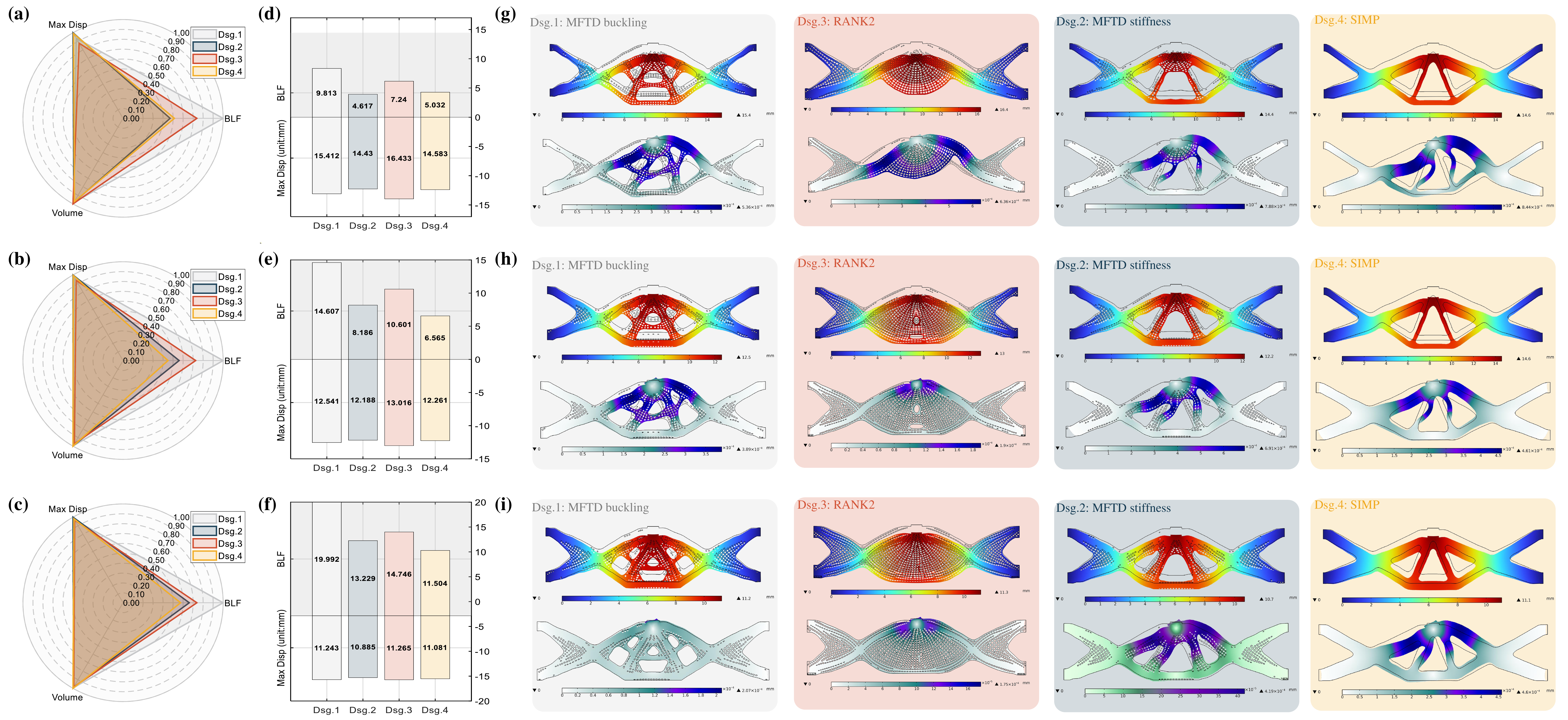}
  \caption{The comparison of the double clamped beam performance for different design strategies. The normalized comprehensive performance evaluation (a)--(c); the detail comparsion for the structural performance (d)--(f); and the detail results for each structure (g)--(i)}\label{fig:The comparison of db}
\end{figure}

To better demonstrate the effectiveness of the proposed method, this section compares the stiffness-based results from Subsection 4.1.2 (Dsg.1: MFTD-stiffness) and the buckling-based results from Subsection 4.1.3 (Dsg.2: MFTD-buckling) with two additional benchmarks: The geometry-mapped results optimized using RANK-2 interpolation (Dsg.3: RANK2), and the traditional SIMP-based solid material distribution optimization (Dsg.4: SIMP). The RANK2 results were reproduced using the open-source code provided in~\cite{woldseth2024808}, while the SIMP results were generated using another open-source implementation from~\cite{andreassen2011efficient}. Except for the low-fidelity optimization process, all other computational settings remain consistent with those used in this study. Figure~\ref{fig:The comparison of db} presents a comparative analysis of the different design approaches (Dsg.1, Dsg.2, Dsg.3, and Dsg.4), focusing on key evaluation metrics, including BLF, volume, and maximum displacement. We selected three sets of results for comparison, ensuring that each set contains designs with similar volumes.

Figure~\ref{fig:The comparison of db}(a), (b), and (c) present radar charts illustrating the normalized performance. The performance distribution can be analyzed by comparing the triangular shape of each radar plot, where a larger coverage area indicates better overall performance. Observations reveal that across all volume levels, Dsg.2 consistently achieves the largest coverage area, followed by Dsg.1 and Dsg.3, with Dsg.4 performing the weakest. Additionally, as volume increases, the coverage areas of Dsg.1 and Dsg.2 remain relatively stable, while those of Dsg.3 and Dsg.4 gradually expand, eventually becoming comparable to Dsg.1.

Figure~\ref{fig:The comparison of db}(g), (h), and (i) present the deformation details of each design method across three different volume levels. In terms of force-induced deformation, the differences between the designs are minimal. Regarding buckling performance, the buckling mode shapes of Dsg.1 and Dsg.4 remain largely consistent, whereas those of Dsg.2 and Dsg.3 vary with changes in volume. A notable observation is that when the volume is around 0.21, both Dsg.2 and Dsg.3 exhibit global buckling, and their buckling areas are nearly identical. However, when the volume increases to around 0.25, the buckling region of Dsg.3 significantly decreases, whereas the buckling region of Dsg.2 remains unchanged. As the volume further increases to around 0.29, Dsg.3 transitions to local buckling, while Dsg.2's buckling region remains almost the same as at 0.25. The detail corresponding numerical values for each metric are presented in Figure~\ref{fig:The comparison of db}(d), (e), and (f).

\subsection{L-bracket beam design} \label{sec:4.2}

The second example involves the design of an L-bracket beam, in which the influence of the mutation operation is investigated. The detailed geometric dimensions and boundary conditions are illustrated in Figure~\ref{fig:The LB beam}(a). The design objective is to minimize the maximum von Mises stress within the structure. Three design strategies are considered: 1. Incorporating the mutation operation with only a volume constraint (Case.1); 2. Without the mutation operation, considering only a volume constraint (Case.2); 3. Incorporating the mutation operation with both volume and stiffness constraints (Case.3).

\subsubsection{Preparation of initial design and optimization settings} \label{sec:4.2.1}
The initial solution generation is performed under the same boundary conditions as the final high-fidelity simulation, with the exception of differences in force magnitude and material properties. Detailed settings are illustrated in Figure~\ref{fig:The LB beam}(b). The low-fidelity optimization model and material parameters follow the same formulation described in Subsection~\ref{sec:4.1.1}, with the objective function still defined as stiffness maximization. A total of 100 samples with varying volume fractions (from 0.2 to 0.4) are generated as initial designs, among which 35 representative examples are selected and presented in Figure~\ref{fig:The LB beam}(c) for reference.

Then, the mutation operation is incorporated into each iteration of the MFTD process. Specifically, ten samples are randomly selected, and for each sample, ten perturbation points \( (i_0, j_0) \) are randomly chosen. The mutation is applied using a spatially decaying influence field characterized by a radius \( r = 20 \) and a smoothness decay factor \( p_0 = 2 \).

\subsubsection{Maximum stress minimization} \label{sec:4.2.2}
\begin{figure}[!tbp]
  \centering
  \includegraphics[width=\textwidth]{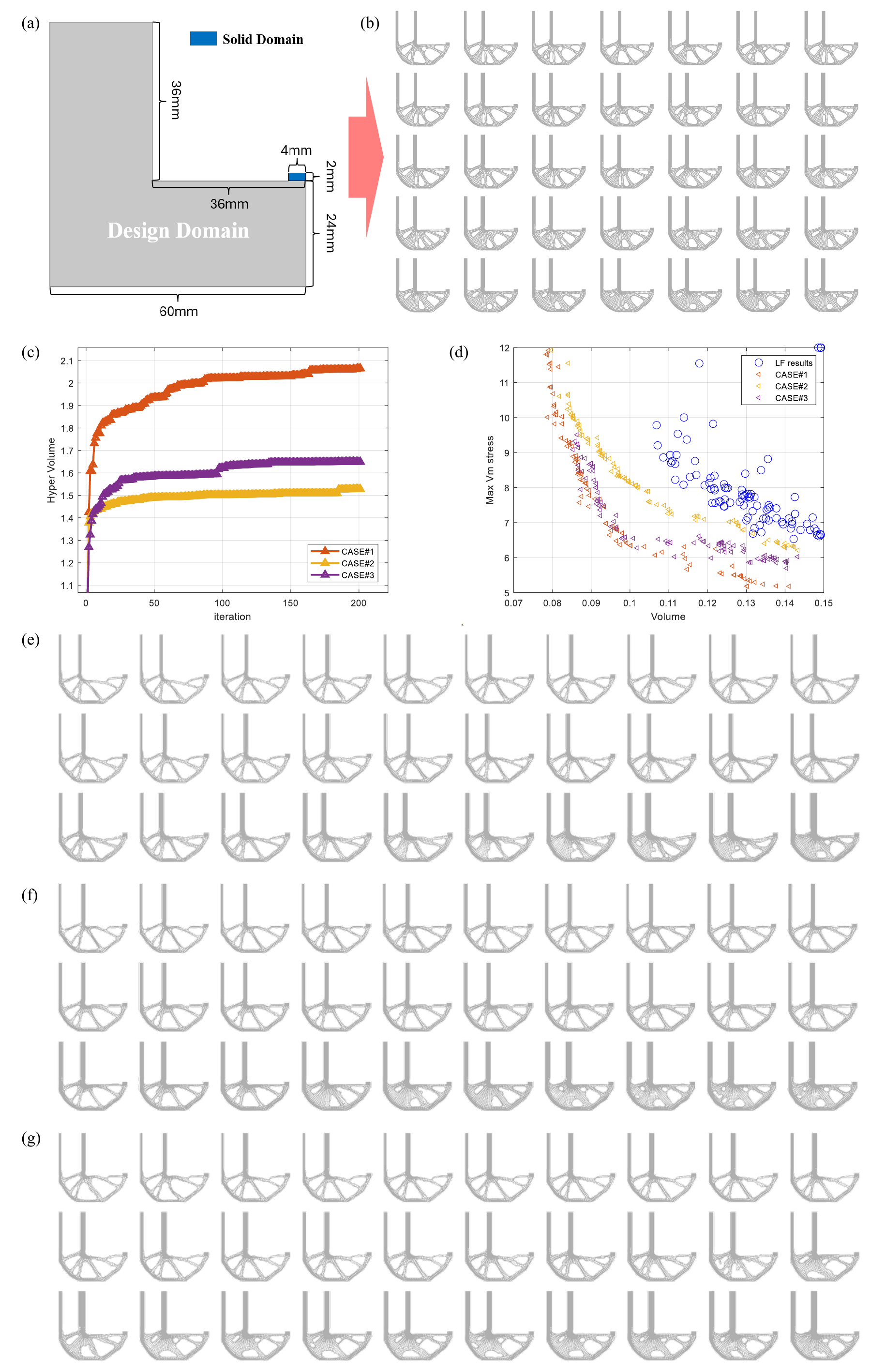}
  \caption{The illustration of design domain and initial designs for L-bracket beam (a)--(b), and the optimized results for three different design strategies (c)--(g).}
  \label{fig:The LB beam}
\end{figure}

Based on the 100 initial designs shown in Figure~\ref{fig:The LB beam}(b), we conducted three separate optimization strategies. Each optimization process consisted of 200 iterations, and the corresponding hypervolume convergence curves are presented in Figure~\ref{fig:The LB beam}(c). Compared to the initial designs, the final optimized hypervolume values exhibit improvements of 107.3\%, 55.1\%, 81.5\% for Case 1, Case 2, and Case 3, respectively. The Pareto fronts of maximum von Mises stress versus volume for the three cases are plotted in Figure~\ref{fig:The LB beam}(d). Subsequently, for each design case, 30 representative MFTD-optimized structures were selected based on their volume distribution and are shown in Figure~\ref{fig:The LB beam}(e), (f), and (g), respectively.

In Case 1, which includes mutation and only considers volume constraints, the final results exhibit noticeably smoother shapes at the corner regions compared to the initial designs, indicating effective stress redistribution. In contrast, Case 2 produces final designs whose corner smoothness is largely consistent with that of the initial configurations. For Case 3, which includes both mutation and a stiffness constraint, the optimized shapes show smoother corners in low-volume cases, whereas high-volume cases retain corner geometries similar to their corresponding initial designs due to the stiffness-preserving constraint.

\subsubsection{Comprehensive analysis for three design strategies} \label{sec:4.2.3}

\begin{figure}[!tbp]
  \centering
  \includegraphics[width=\textwidth]{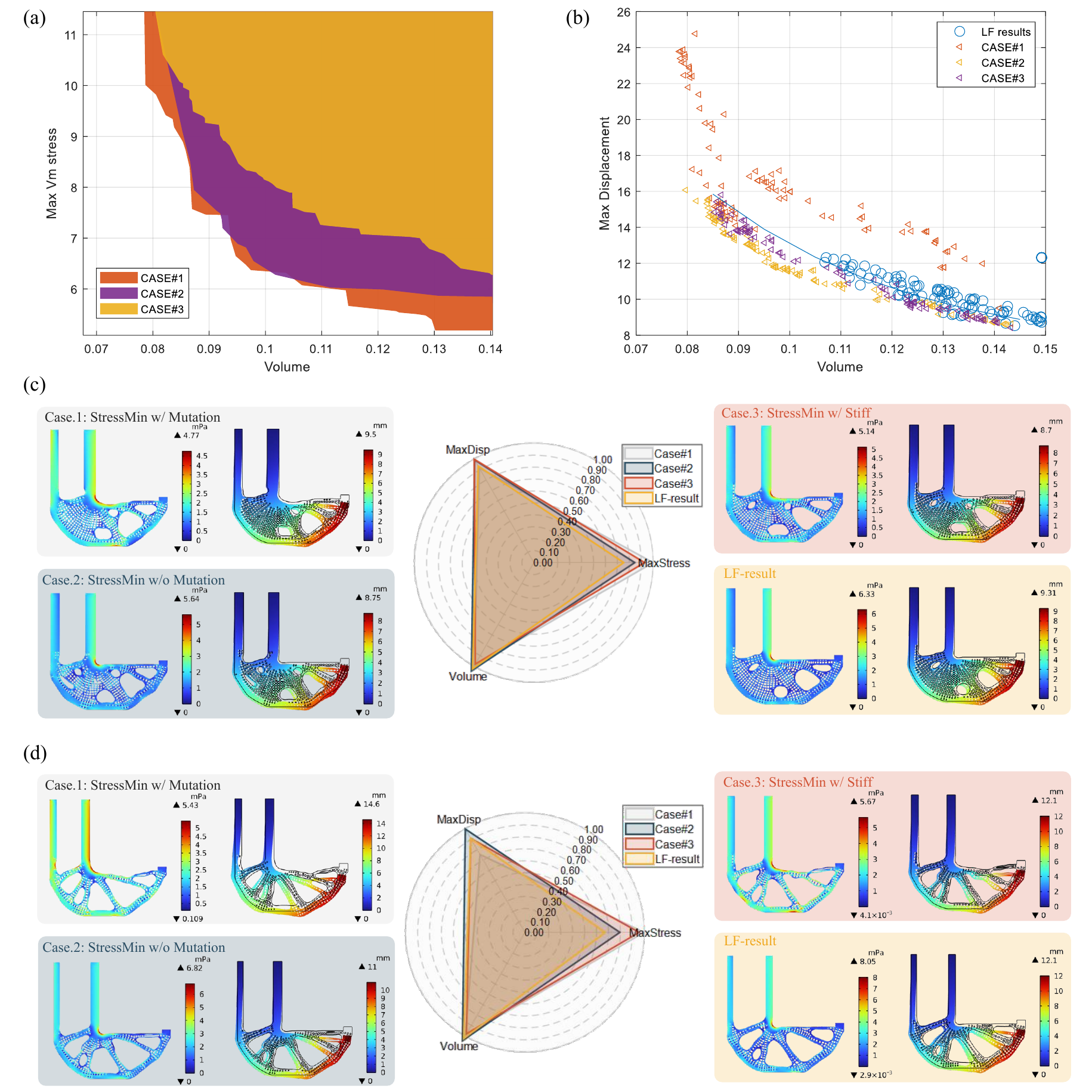}
  \caption{The comprehensive comparsion for the L-bracket beam results with different design strategies.}
  \label{fig:The LB beam 2}
\end{figure}

In this subsection, we perform a comprehensive comparison of structural performance across the three design strategies (Case 1, Case 2, and Case 3). As shown in Figure~\ref{fig:The LB beam 2}(a), we illustrate the distribution of maximum von Mises stress versus volume for all cases. The visualization is based on hypervolume dominance coverage, representing the performance envelope of each design case. It is evident that Case 1 achieves the largest coverage area, followed by Case 2, while Case 3 exhibits the smallest coverage. This indicates that Case 1, which incorporates volume constraints along with mutation operations, demonstrates the best performance in terms of minimizing stress under varying volume fractions.

Subsequently, Figure~\ref{fig:The LB beam 2}(b) presents the distribution of maximum displacement versus volume, serving as an indicator of structural stiffness. The blue stiffness line denotes the average stiffness of the initial solutions and acts as a reference lower bound. From the figure, it can be observed that for Case 1, which only imposes volume constraints, most designs exhibit a noticeable degradation in stiffness compared to the initial solutions. In contrast, Case 2 (without mutation) and Case 3 (which includes both mutation and explicit stiffness constraints) result in solutions that largely comply with the stiffness requirement, maintaining displacements close to or below the reference bound.

Finally, two groups of results were selected for comparative analysis across four structural designs. Within each group, the structural volumes were approximately equal. Figure~\ref{fig:The LB beam 2}(c) and (d) present radar charts illustrating the normalized performance metrics, including maximum stress, stiffness (evaluated via displacement), and volume. The overall performance distribution can be interpreted by examining the shape and coverage area of each radar plot—larger enclosed areas indicate superior overall performance across the three metrics.

According to the observations, when the structural volume is relatively small, there is a noticeable difference among the four design strategies. In particular, case.3, which incorporates mutation operations and includes a stiffness constraint, consistently demonstrates superior overall performance across stress, stiffness, and volume. As the volume increases, the distinctions among the four strategies become less pronounced. Nevertheless, across all volume levels, the LF-based results consistently exhibit the poorest performance.

\subsection{Multiple objective design} \label{sec:4.3}

For non-gradient optimization, multi-objective optimization has always been one of its strengths. In this case, we consider the design of a typical part under multiple loading conditions. The boundary conditions and structural dimensions are illustrated in Figure~\ref{fig:The initial design of part}(a). 

\subsubsection{Preparation of initial design} \label{sec:4.3.1}

In this case study, we aim to optimize the internal infill structure while preserving the overall external shape of the part to suit multi-loading condition designs. Therefore, in the low-fidelity optimization model, a new design variable \( \rho \in [0,1] \) is introduced to independently interpolate the width of the structural members:

\begin{equation}
  \mu = \mu_{min}+(\mu_{max}-\mu_{min})\rho^p.
\end{equation}

The form of this design variable follows the traditional SIMP interpolation scheme, where \(\mu_{\text{min}}\) represents the minimum width of the structural members across the design domain, and \(\mu_{\text{max}}\) denotes the maximum width. In this study, \(\mu_{\text{min}} = 0.3\) and \(\mu_{\text{max}} = 1\). This interpolation is defined such that when the design variable (or density) \(\rho = 0\), the member width reaches its minimum value, corresponding to the thinnest structural members; and when \(\rho = 1\), it corresponds to a fully solid region. Accordingly, the effective elasticity tensor is modified to incorporate the influence of \( \rho \) as

\begin{equation}
\tilde{\mathbf{C}}^{\text{H}}(\theta, \rho) = \mathbf{T}^\top(\theta) \, \mathbf{C}^{\text{H}}(\mu(\rho), \mu(\rho)) \, \mathbf{T}(\theta),
\label{eq:macro-rho}
\end{equation}

The geometric variable vector is now defined as
\begin{equation}
  \mathbf{X} = \begin{bmatrix} \mu(\rho) & \mu(\rho) & \theta \end{bmatrix}^\top.
\end{equation}

It involves only two design variables: \(\rho\) and \(\theta\). Following Eq.\eqref{eq:low-fidelity}, the specified optimization model is formulated as

\begin{equation}
  \begin{array}{ll}
  \text{Find:} &  \theta, \rho \\
  \text{Minimize:} & \text{obj}(\theta, \rho) \\
  \text{Subject to:} & \mathbf{K}(\theta, \rho) \mathbf{U} = \mathbf{F} \\
  & V(\rho) \leq V_0 \\
  & 0 \leq \rho \leq 1 \\
  & -4\pi \leq \theta \leq 4\pi
  \end{array}
\end{equation}

\begin{figure}[!tbp]
  \centering
  \includegraphics[width=\textwidth]{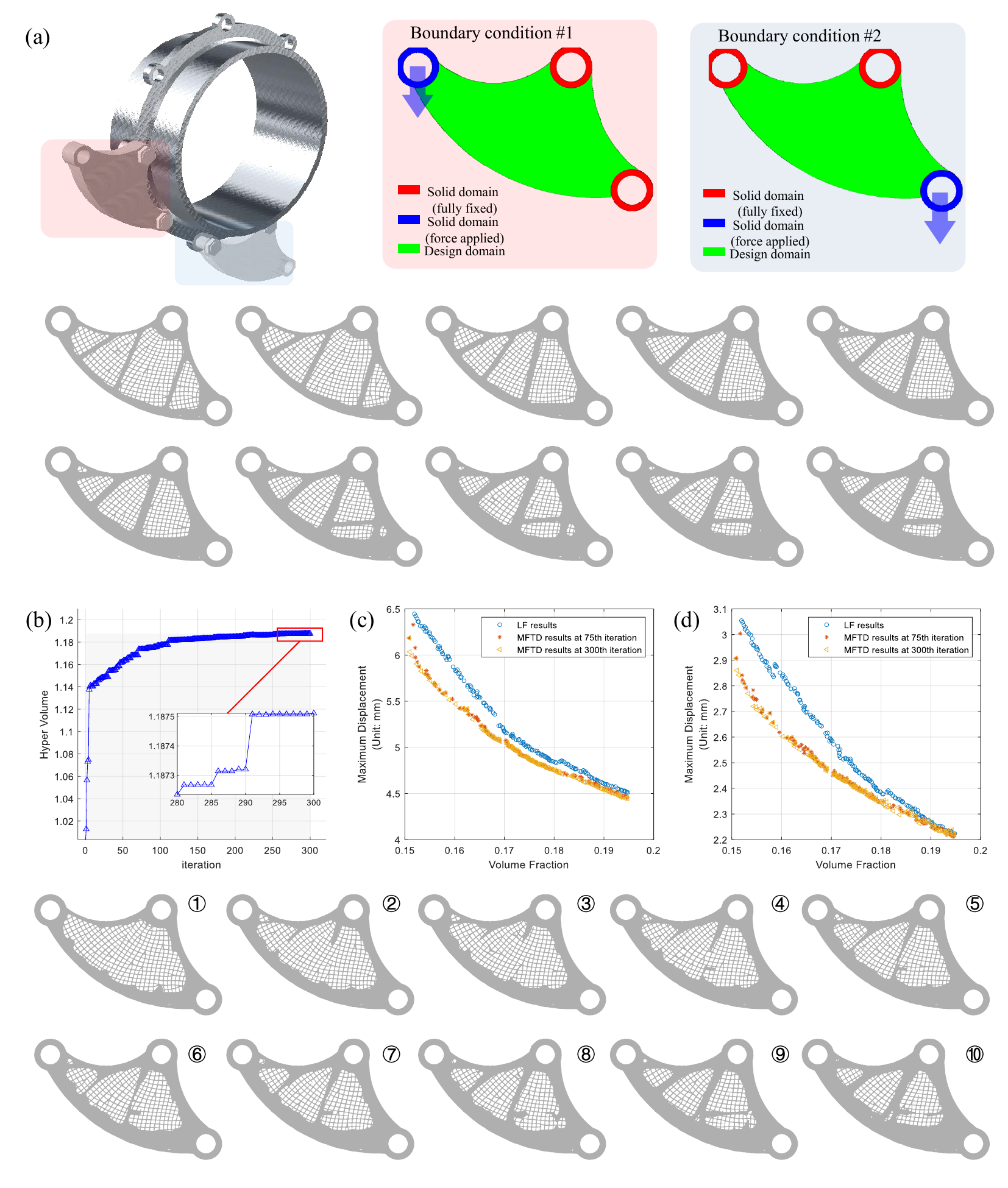}
  \caption{The optimization result for the part design. (a) is the design domain and boundary conditions, and (b) is the hypervolume value convergence history, and (c)(d) are the corresponding pareto front of the responses under different boundary conditions. 10 selected optimized designs are shown below them..}\label{fig:The initial design of part}
\end{figure}

The initial solution generation follows the same boundary conditions as the final high-fidelity simulation, except for differences in force magnitude and material properties. Their information is the same with double clamper beam case. Detailed information is provided in Figure~\ref{fig:The initial design of part}(b). In this case, only 160 initial designs are created according to the change of $V_0$(0.3 to 0.6). In Figure~\ref{fig:The initial design of part}, we present 10 representative initial solutions, and their CAD model after geometry mapping.

\subsubsection{Stiffness improvement} \label{sec:4.3.2}

Based on the initial solutions shown in Subsection 4.1.1, we optimize the structure with the objective of minimizing maximum displacement under different boundary conditions and volume. The optimization process still consists of 300 iterations, and the corresponding hypervolume convergence curve is shown in Figure~\ref{fig:The initial design of part}(b). Compared to the initial solution, the final optimized hypervolume shows a 18.75\% improvement. The Pareto front of Maximum Displacement vs. Volume for two boundary conditions are plotted in Figure~\ref{fig:The initial design of part}(c) and (d). 

According to observations from Figure~\ref{fig:The initial design of part}(c) and (d), it can be seen that under both boundary conditions, as the volume increases, the maximum displacement gradually decreases, demonstrating a typical optimization trade-off. When comparing at the same volume, significant improvements are observed with the MFTD results at lower volume fractions. However, as the volume fraction increases, the extent of improvement diminishes. For Boundary Condition 1, a slight improvement is still maintained, whereas for Boundary Condition 2, the MFTD results become nearly identical to the initial solution.

10 representative MFTD-optimized structures are selected to show in Figure~\ref{fig:The initial design of part} based on their volume distribution. Based on observations, it can be seen that when the volume fraction is very low, the design tends to favor pure lattice infill structures. As the volume fraction increases, an increasing amount of solid regions appears. This indicates that under relatively complex loading conditions, lattice infill designs may be more efficient than solid designs at low volume fractions.

\subsubsection{Comprehensive analysis for the design results} \label{sec:4.3.3}

\begin{figure}[!tbp]
  \centering
  \includegraphics[width=\textwidth]{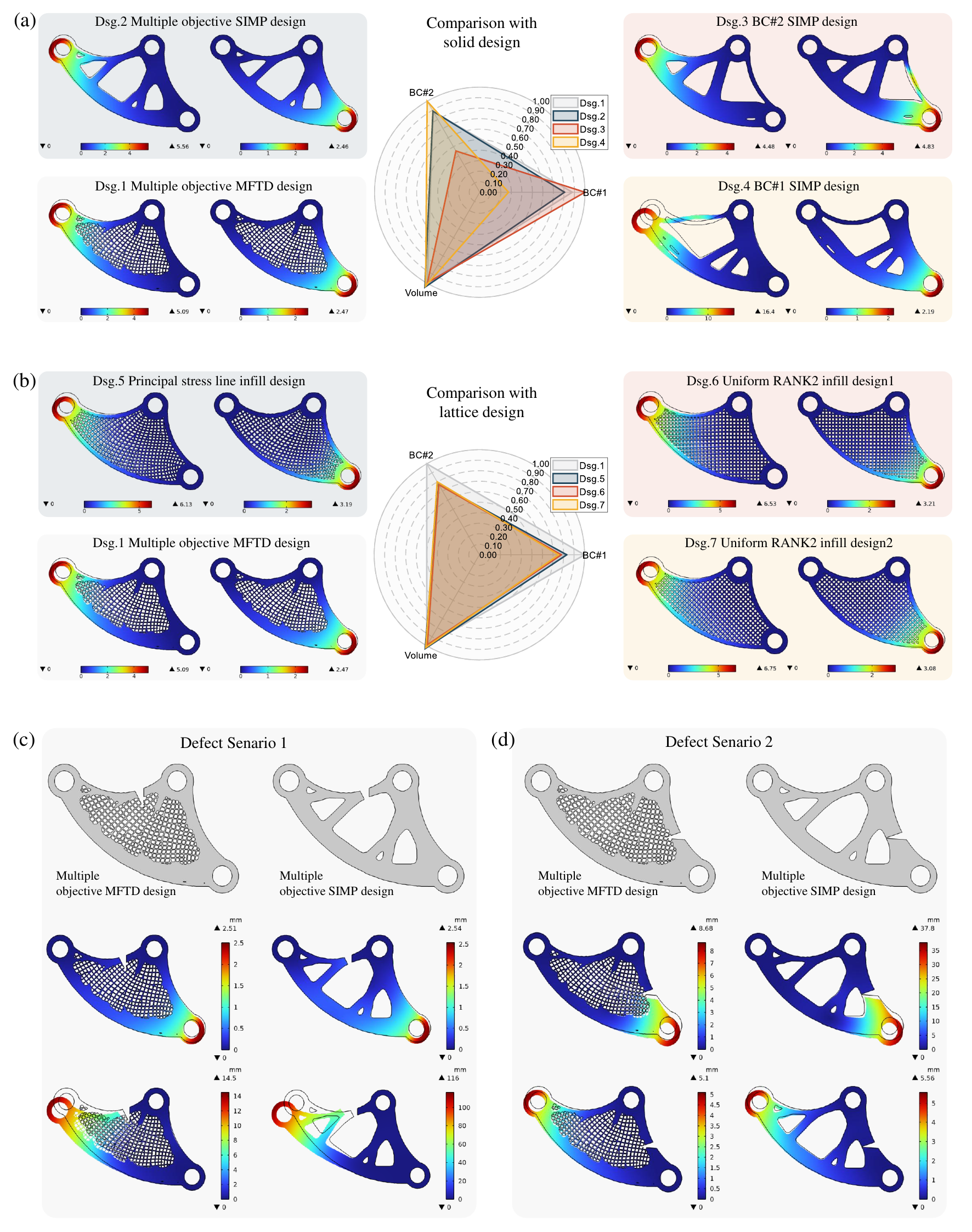}
  \caption{The comparison of structure performance for different design. (a) The performance comparison within the solid design group; (b) the performance comparison within the lattice design group; and (c)--(d) the performance comparison under different defect scenarios}\label{fig:The result of part}
\end{figure}

To better demonstrate the effectiveness of the proposed method, we firstly designed two sets of comparative experiments: one based on solid designs and the other based on lattice designs. The MFTD results were compared against these benchmarks while maintaining approximately the same volume. For the solid design control group, three different cases were considered: a multi-objective SIMP result optimized for both Boundary Condition 1 and Boundary Condition 2, a SIMP result optimized solely for Boundary Condition 1, and a SIMP result optimized solely for Boundary Condition 2. For the lattice design control group, we selected three representative lattice infill strategies: RANK2 filling based on principal stress directions, RANK2 filling with material orientation at \(0^\circ\), and RANK2 filling with material orientation at \(45^\circ\).

The corresponding analysis results are shown in Figure~\ref{fig:The result of part}(a) and (b). For the solid design comparison group, it can be observed that each design (Dsg.3 and Dsg.4) achieves superior performance under its respective optimized boundary condition, outperforming both the MFTD results (Dsg.1) and the multi-objective SIMP results (Dsg.2). However, under non-optimized loading conditions, their performance degrades significantly, highlighting the lack of generalization capability inherent in single-objective SIMP-based designs. When directly comparing the MFTD results with the multi-objective SIMP results, it is found that under Boundary Condition 2, both approaches achieve largely comparable structural performance. Nevertheless, under Boundary Condition 1, the MFTD method exhibits a clear advantage, demonstrating better stress distribution and lower peak stresses. Radar charts summarizing the normalized performance of each design method indicate that the multi-objective SIMP designs cover a broader performance area than the single-objective SIMP designs. Moreover, the MFTD designs exhibit an even larger coverage area compared to the multi-objective SIMP results, reflecting stronger robustness and better adaptability across different loading conditions.

In the lattice design comparison group, the performance differences among various designs under different boundary conditions are relatively moderate. Nevertheless, the MFTD design consistently exhibits the best overall performance, balancing structural strength, material efficiency, and robustness across all tested conditions.

Finally, an additional set of experiments was conducted to evaluate the robustness of each optimization strategy under artificial damage scenarios. In these tests, the generated structures were deliberately subjected to severe local damages (as illustrated in Figure~\ref{fig:The result of part}(c) and (d)) to simulate critical failure conditions. Finite element analyses were performed to assess their mechanical performance post-damage. The results clearly show that the MFTD designs maintain superior mechanical integrity under both types of damage compared to the multi-objective SIMP designs, further highlighting the enhanced robustness achieved by the MFTD framework.

\subsection{Computational cost analysis} \label{sec:4.4}

In this subsection, we analyze three representative results selected from the case.1 (L-bracket beam) shown in Subsection~\ref{sec:4.1.2}, all of which share the same boundary conditions and parameter settings but differ in their decision variable counts. Three cases are considered: 100 samples, 200 samples, and 300 samples.

\begin{figure}[!tbp]
  \centering
  \centering
  \includegraphics[width=\textwidth]{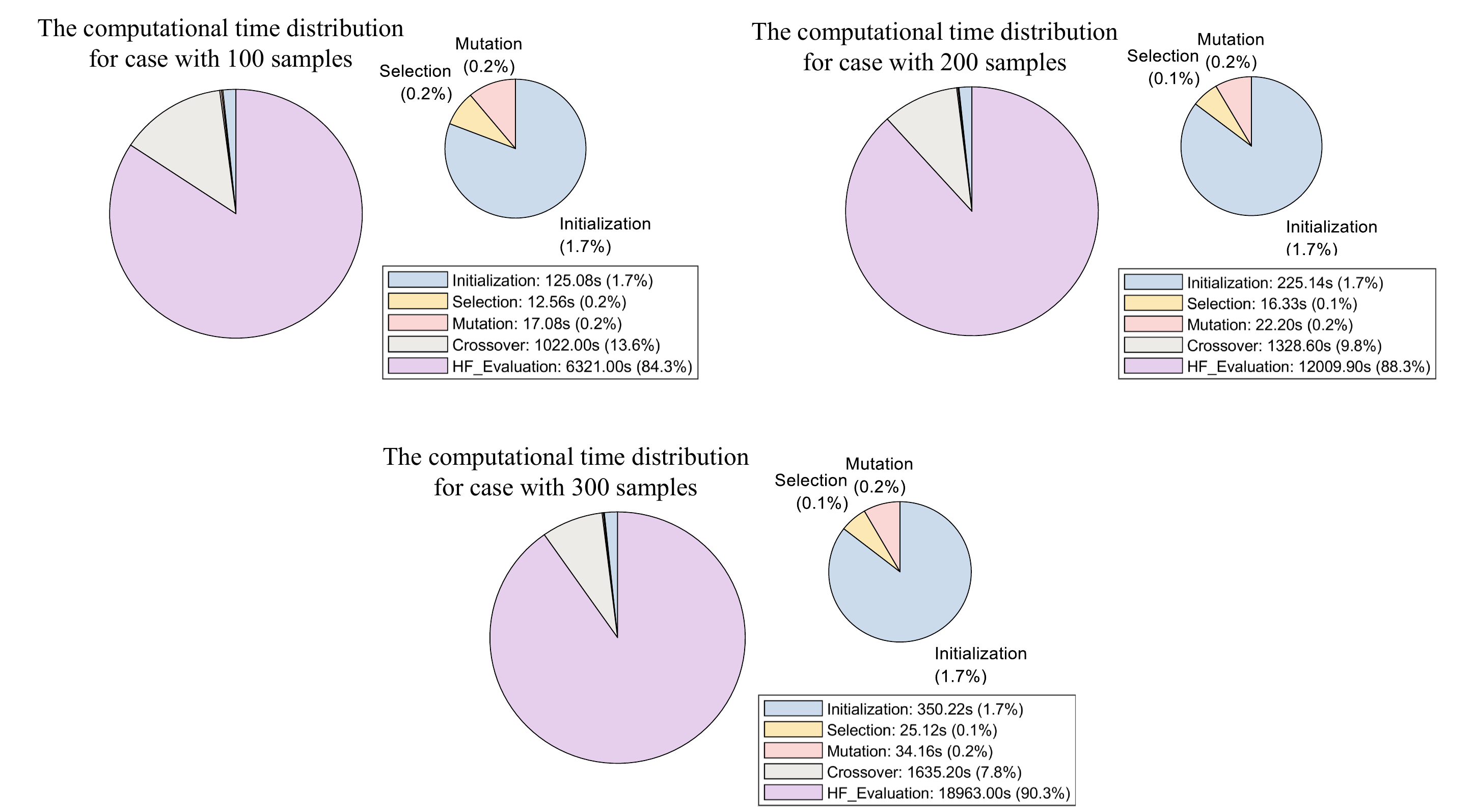} 
  \caption{Summary of Computational Time.}
  \label{fig:times}
\end{figure}

The three pie charts in Figure~\ref{fig:times} summarize the total computational time for each case, as well as the distribution of time spent across different components. The computation time is categorized according to the modules defined in Subsection~\ref{sec:Optimization model}, which include: (1) obtaining the initial solution through low-fidelity optimization, (2) generating and analyzing the high-fidelity model, (3) selecting solutions using the NSGA-II algorithm, (4) using image processing to achieve mutation operation, and (5) using PCA to compress the data and then training the VAE and generating new candidate solutions.

To ensure consistency in comparison, we fixed the number of optimization iterations to 100 for all cases. As observed, the overall computational cost increases with the number of samples: Case with 100 samples took 7497.7 seconds, Case with 200 samples took 13602 seconds, and Case with 300 samples took 21008 seconds. However, the increase in model training cost is relatively modest, rising only from 1022 seconds to 1635 seconds. In fact, its proportional share of the total time decreased from 13.6\% to 7.8\%. In contrast, the computational cost for low-fidelity data generation and candidate selection is negligible and can be largely ignored. The high-fidelity analysis, on the other hand, grows significantly with sample size and constitutes the primary computational bottleneck.

\section{Concluison} \label{sec:Concluison}

{This study proposed a design framework, which seamlessly integrates MFTD with phasor-based de-homogenization techniques, and enables efficient optimization of complex multiscale structures with different function considerations.} Within the proposed framework, PCA was used to reduce the dimensionality of the low-fidelity design variables, while a phasor-based de-homogenization method was employed to generate detailed manufacturable geometries. A data-driven evolutionary algorithm, incorporating NSGA-II for selection, a VAE for crossover operations, and a image processing method for mutation, was developed to efficiently search the latent design space. This approach successfully addresses the inherent challenges of non-differentiability and nonlinearity in the mapping between abstract design variables and full-resolution geometries.

Through a series of numerical case studies, including different design objectives: stiffness (Subsection~\ref{sec:4.1.2}), buckling (Subsection~\ref{sec:4.1.3}), stress (Subsection~\ref{sec:4.2.2}), and multi-loading condition design (Subsection~\ref{sec:4.3.2}); and different boundary conditions: double clamper beam (Subsection~\ref{sec:4.1}), L-bracket beam (Subsection~\ref{sec:4.2}), and random part (Subsection~\ref{sec:4.3}), the effectiveness and robustness of the proposed method were validated. 

Compared to traditional SIMP-based and homogenization-based lattice designs, the evolutionary de-homogenization method achieved superior performance in terms of multi-objective trade-offs, geometric fidelity, and structural robustness, particularly under complex loading conditions and damage scenarios. Additionally, the computational cost analysis demonstrated that although high-fidelity simulations constitute the primary bottleneck, the proposed method's generative model training cost remains relatively low due to the reduced variable dimension, further ensuring practical feasibility for large-scale design problems.

Despite the promising results, there are several avenues for future work. One important direction is the experimental validation of the fabricated structures. Moreover, extending the current two-dimensional framework to fully three-dimensional cases would be essential for broader engineering applications. Further improvements in generative modeling efficiency, such as adopting more advanced generative model, could also enhance scalability and generalization.

\bibliographystyle{unsrt} 
\bibliography{Reference} 

\end{document}